\documentclass[preprint]{elsarticle}
\usepackage{hyperref}
\usepackage{latexsym}
\usepackage{graphicx}
\usepackage{multicol,multirow}
\usepackage{amsmath,amssymb,amsfonts}
\usepackage{mathrsfs}
\usepackage{amsthm}
\usepackage{subcaption}
\usepackage{upgreek}
\usepackage{xcolor}

\newtheorem*{remark}{Remark}

\usepackage{xcolor}

\definecolor{Blue}{rgb}{0,0,1}
\definecolor{Red}{rgb}{1,0,0}
\newcommand{\XH}[1]{{\textcolor{Red}{#1}}}

\renewcommand{\textcolor}[2]{#2} 
\journal{Elsevier}

\usepackage{numcompress}

\begin{document}

\begin{frontmatter}

\title{Thermodynamically Consistent Machine-Learned Internal State Variable Approach for Data-Driven Modeling of Path-Dependent Materials}



\author[1]{Xiaolong He}
\author[1]{Jiun-Shyan Chen\corref{corresponding_author}}
\cortext[corresponding_author]{Corresponding author}
\ead{js-chen@ucsd.edu}

\address[1]{Department of Structural Engineering, University of California, San Diego, La Jolla, CA, 92093, USA}





\begin{abstract}
Characterization and modeling of path-dependent behaviors of complex materials by phenomenological models remains challenging due to difficulties in formulating mathematical expressions and internal state variables (ISVs) governing path-dependent behaviors. Data-driven machine learning models, such as deep neural networks and recurrent neural networks (RNNs), have become viable alternatives. 
However, pure black-box data-driven models mapping inputs to outputs without considering the underlying physics suffer from unstable and inaccurate generalization performance. 
This study proposes a machine-learned physics-informed data-driven constitutive modeling approach for path-dependent materials based on the measurable material states. 
The proposed data-driven constitutive model is designed with the consideration of universal thermodynamics principles, where the ISVs essential to the material path-dependency are inferred automatically from the hidden state of RNNs. The RNN describing the evolution of the data-driven machine-learned ISVs follows the thermodynamics second law.
To enhance the robustness and accuracy of RNN models, stochasticity is introduced to model training.
The effects of the number of RNN history steps, the internal state dimension, the model complexity, and the strain increment on model performances have been investigated.
The effectiveness of the proposed method is evaluated by modeling soil material behaviors under cyclic shear loading using experimental stress-strain data.
\end{abstract}
\end{frontmatter}


\section{Introduction}\label{sec:intro}
Traditional constitutive modeling is based on constitutive or material laws to describe the explicit relationship among the measurable material states, e.g., stresses and strains, and internal state variables (ISVs) based on experimental observations, mechanistic hypothesis, and mathematical simplifications.
However, limited data and functional form assumptions inevitably introduce errors to the model parameter calibration and model prediction.
Moreover, with the pre-defined functions, constitutive laws often lack generality to capture full aspects of material behaviors \cite{HeJBM2020,he2021deep}.

Path-dependent constitutive modeling typically applies models with evolving ISVs in addition to the state space of deformation \cite{coleman1967thermodynamics,horstemeyer2010historical}. 
The ISV constitutive modeling framework has been effectively applied to model various nonlinear solid material behaviors, e.g., elasto-plasticity \cite{kratochvil1969thermodynamics,simo1992associative}, visco-plasticity \cite{simo1992algorithms}, and material damage \cite{perzyna1986internal}.
However, ISVs are often non-measurable, which makes it challenging to define a complete and appropriate set of ISVs for highly nonlinear and complicated materials, e.g., geomechanical materials. Further, the traditional ISV constitutive modeling approach often results in excessive complexities with high computational cost, which is undesirable in practical applications.

\XH{In recent years, machine learning (ML) based data-driven approaches have demonstrated successful applications in various engineering problems, such as solving partial differential equations \cite{raissi2019physics,he2021physics,karniadakis2021physics,kadeethum2021framework}, 
system or parameter identification \cite{brunton2016discovering,raissi2019physics,cranmer2020discovering,tartakovsky2020physics,haghighat2021physics,kadeethum2021framework,he2022physics}, 
data-driven computational mechanics \cite{kirchdoerfer2017data,ayensa2018new,he2019physics,eggersmann2019model,he2020physics,he2021deep,kanno2021kernel,bahmani2021manifold}, 
reduced-order modeling \cite{xie2019non,bai2021non,kaneko2021hyper,kim2022fast,fries2022lasdi,he2022glasdi,kadeethum2022non}, 
material design \cite{bessa2017framework,butler2018machine}, etc.
ML models, such as deep neural networks (DNNs), have emerged as a promising alternative for constitutive modeling due to their strong flexibility and capability in extracting complex features and patterns from data \cite{bock2019review}}.
DNNs have been applied to model a variety of materials, including concrete materials \cite{ghaboussi1991knowledge}, hyper-elastic materials \cite{Shen2005}, visco-plastic material of steel \cite{furukawa1998implicit}, and homogenized properties of composite structures \cite{lefik2009artificial}. DNN-based constitutive models haven been integrated into finite element solvers to predict path- or rate-dependent materials behaviors \cite{lefik2003artificial,hashash2004numerical,jung2006neural,Stoffel2019,zhang2020using}. Recently, physical constraints or principles have been integrated into DNNs for data-driven constitutive modeling, including symmetric positive definiteness \cite{xu2021learning}, material frame invariance \cite{ling2016machine}, and thermodynamics \cite{vlassis2021sobolev,masi2021thermodynamics}. However, to model path-dependent materials, the DNN-based constitutive models require fully understood and prescribed material's internal states, which is difficult for materials with highly nonlinear and complicated path-dependent behaviors and limits their applications in practice.

Recurrent neural networks (RNNs) designed for sequence learning have been successfully applied in various domains, such as machine translation and speech recognition, due to their capability of learning history-dependent features that are essential for sequential prediction \cite{lipton2015critical,yu2019review}. The RNN and gated variants, e.g., the long short-term memory (LSTM) \cite{hochreiter1997long} cells and the gated recurrent units (GRUs) \cite{cho2014properties, chung2014empirical}, have been applied to path-dependent materials modeling \cite{heider2020so}, including plastic composites \cite{mozaffar2019deep}, visco-elasticity \cite{chen2021recurrent}, and homogeneous anisotropic hardening \cite{gorji2020potential}. RNN-based constitutive models have also been applied to accelerate multi-scale simulations with path-dependent characteristics \cite{wang2018multiscale,ghavamian2019accelerating,wu2020recurrent,logarzo2021smart,wu2022recurrent}. Recently, Bonatti and Mohr \cite{bonatti2022importance} proposed a self-consistent RNN for path-dependent materials such that the model predictions converge as the loading increment is decreased.
However, these RNN-based data-driven constitutive models may not satisfy the underlying thermodynamics principles of path-dependent materials.

In this study, we propose a thermodynamically consistent machine-learned ISV approach for data-driven modeling of path-dependent materials, which relies purely on measurable material states. The first thermodynamics principle is integrated into the model architecture whereas the second thermodynamics principle is enforced by a constraint on the network parameters. In the proposed model, an RNN is trained to infer intrinsic ISVs from its hidden (or memory) state that captures essential history-dependent features of data through a sequential input. The RNN describing the evolution of the data-driven machine-learned ISVs follows the thermodynamics second law. In addition, a DNN is trained simultaneously to predict the material energy potential given strain, ISVs, and temperature (for non-isothermal processes). Further, model robustness and accuracy is enhanced by introducing \textit{stochasticity} to inputs for model training to account for uncertainties of input conditions in testing.

The remainder of this paper is organized as follows. The background of thermodynamics principles is first introduced in Section \ref{sec:thermodynamics}.
In Section \ref{sec:rnn}, DNNs and RNNs are introduced and their applications to path-dependent materials modeling are discussed. 
Section \ref{sec:tcrnn} introduces the proposed thermodynamically consistent machine-learned ISV approach for data-driven modeling of path-dependent materials, where two thermodynamically consistent recurrent neural networks (TCRNNs) are discussed.
Finally, in Section \ref{sec:result}, the effectiveness and generalization capability of the proposed TCRNN models are examined by modeling an elasto-plastic material and undrained soil under cyclic shear loading. A parametric study is conducted to investigate the effects of the number of RNN steps, the internal state dimension, the model complexity, and the strain increment on the model performance.
Concluding remarks and discussions are summarized in Section \ref{sec:conclusion}.
\section{Thermodynamics Principles}\label{sec:thermodynamics}

The balance of energy, i.e., the \textit{first thermodynamics principle}, can be expressed as \cite{silhavy2013mechanics}
\begin{equation}\label{eq.1st_law}
    \rho \dot{e} = \boldsymbol{\sigma} : \dot{\boldsymbol{\varepsilon}} - \text{div } \mathbf{q} + \rho h,
\end{equation}
where \XH{the superposed ``." denotes the material time derivative;} $\rho$ is the material density; $e$ is the specific internal energy and $\dot{e}$ is its rate; $\boldsymbol{\sigma}$ is the Cauchy stress tensor; $\boldsymbol{\varepsilon}$ is the strain tensor; $\boldsymbol{\sigma} : \dot{\boldsymbol{\varepsilon}}$ denotes the rate of mechanical work; $\mathbf{q}$ is the heat flux; $h$ is the specific rate of heat supply. The local form of the \textit{second thermodynamics principle} expressed in terms of the Clausius-Duhem inequality reads \cite{silhavy2013mechanics}
\begin{equation}\label{eq.2nd_law}
    \rho\dot{s} + \text{div } \big( \frac{\mathbf{q}}{T} \big) - \frac{\rho h}{T} \ge 0,
\end{equation}
where $s$ denotes the specific entropy and $T$ is the positive absolute temperature. 
Combining the first and second thermodynamic principles yields the dissipation inequality
\begin{equation}\label{eq.total_disspation_inequality}
    \boldsymbol{\sigma} : \dot{\boldsymbol{\varepsilon}} - \rho \dot{e} + T \rho \dot{s} - \mathbf{q} \cdot \frac{\nabla T}{T} \ge 0.
\end{equation}
The left-hand side of Eq. (\ref{eq.total_disspation_inequality}) represents the total dissipation rate that can be decomposed into the \XH{non-negative} mechanical dissipation rate $D$ and the \XH{non-negative} thermal dissipation rate $D^{th}$ \cite{silhavy2013mechanics,simo2006computational}: 
\begin{subequations}\label{eq.mech_thermal_disspation_rate}
    \begin{align}
        & D = \boldsymbol{\sigma} : \dot{\boldsymbol{\varepsilon}} - \rho \dot{e} + T \rho \dot{s} \ge 0, \label{eq.mech_diss_rate}\\
        & D^{th} = - \mathbf{q} \cdot \frac{\nabla T}{T} \ge 0. \label{eq.therm_diss_rate}
    \end{align}
\end{subequations}
Considering a constant material density and defining the specific internal energy per unit volume as $E = \rho e$ and the specific entropy per unit volume as $S = \rho s$, we have $\dot{E} = \rho \dot{e}$ and $\dot{S} = \rho \dot{s}$. Therefore, Eq. (\ref{eq.mech_diss_rate}) can be rewritten as
\begin{equation}\label{eq.mech_disspation_rate_inequality2}
    D = \boldsymbol{\sigma} : \dot{\boldsymbol{\varepsilon}} - \dot{E} + T \dot{S} \ge 0.
\end{equation}

Denoting the specific Helmholtz free energy as $F = E - TS$ \cite{silhavy2013mechanics} and taking the time derivative gives
\begin{equation}\label{eq.free_energy_rate}
    \dot{F} = \dot{E} - \dot{T}S - T\dot{S}.
\end{equation}
Combining Eq. (\ref{eq.mech_disspation_rate_inequality2}) and Eq. (\ref{eq.free_energy_rate}) gives
\begin{equation}\label{eq.free_energy_dissipation}
    \dot{F} = \boldsymbol{\sigma} : \dot{\boldsymbol{\varepsilon}} - D - \dot{T} S.
\end{equation}
For strain-rate independent materials, the Helmholtz free energy can be defined as \cite{silhavy2013mechanics}
\begin{equation}\label{eq.free_energy_definition}
    F := F(T, \boldsymbol{\varepsilon}, \mathbf{z}),
\end{equation}
where $\mathbf{z} = (z_1, ..., z_N)$ is a collection of $N$ internal state variables (ISVs) introduced to characterize the state of path-dependent materials, which can also be interpreted as history variables \cite{eggersmann2019model}. However, ISVs are often non-measurable and the identification of the ISVs is often based on empiricism, which is non-trivial for materials with highly complex and nonlinear path-dependent behaviors. Here, a ML-enhanced data-driven approach is proposed to automatically infer the essential ISVs that follow the thermodynamics principles, which will be discussed in Section \ref{sec:tcrnn}. Differentiation of Eq. (\ref{eq.free_energy_definition}) gives
\begin{equation}\label{eq.free_energy_diff}
    \dot{F} = \frac{\partial F}{\partial T} \dot{T} + \frac{\partial F}{\partial \boldsymbol{\varepsilon}}: \dot{\boldsymbol{\varepsilon}} + \frac{\partial F}{\partial \mathbf{z}} \cdot \dot{\mathbf{z}}.
\end{equation}
Equating Eq. (\ref{eq.free_energy_dissipation}) with Eq. (\ref{eq.free_energy_diff}) gives
\begin{equation}\label{eq.equilibrium}
    \bigg( \frac{\partial F}{\partial T} + S \bigg) \dot{T} + \bigg( \frac{\partial F}{\partial \boldsymbol{\varepsilon}} - \boldsymbol{\sigma} \bigg) : \dot{\boldsymbol{\varepsilon}} + \bigg( \frac{\partial F}{\partial \mathbf{z}} \cdot \dot{\mathbf{z}} + D \bigg) = 0.
\end{equation}
The arbitrariness of $\dot{T}$, $\dot{\boldsymbol{\varepsilon}}$, and $\dot{z}$ leads to the following relations
\begin{subequations}\label{eq.S_sigma_D}
    \begin{align}
        S & = - \frac{\partial F}{\partial T}, \\
        \boldsymbol{\sigma} & = \frac{\partial F}{\partial \boldsymbol{\varepsilon}}, \\
        D & = - \frac{\partial F}{\partial \mathbf{z}} \cdot \dot{\mathbf{z}}.
    \end{align}
\end{subequations}
These relations are derived based on the universal thermodynamics principals and considered in the proposed data-driven constitutive models. In Section \ref{sec:tcrnn}, we will introduce a thermodynamically consistent machine-learned ISV approach for data-driven modeling of path-dependent materials with the consideration of the thermodynamically consistent relations (Eq. (\ref{eq.S_sigma_D})).
\section{Black-Box Data-Driven Modeling of Path-Dependent Materials}\label{sec:rnn}

\subsection{Deep Neural Networks (DNNs) Constitutive Models}\label{sec:rnn-nn}

Deep neural networks (DNNs), as the core of the deep learning \cite{goodfellow2016deep}, represent complex models that relate data inputs, $\mathbf{x} \in \mathbb{R}^{d_{in}}$, to data outputs, $\mathbf{y} \in \mathbb{R}^{d_{out}}$. 
A typical DNN is composed of an input layer, an output layer, and $L$ hidden layers. Each hidden layer transforms the outputs of the previous layer through an affine mapping followed by a nonlinear activation function $a(\cdot)$, which can be written as:
\begin{equation}\label{eq.nn}
    \mathbf{x}^{(l)} = a(\mathbf{W}^{(l)} \mathbf{x}^{(l-1)} + \mathbf{b}^{(l)}), \hspace{0.5cm} l = 1,...,L,
\end{equation}
where the superscript $(l)$ denotes the layer the quantities belong to, e.g., $\mathbf{x}^{(l)} \in \mathbb{R}^{n_l}$ is the outputs of layer $l$ with $n_l$ neurons.
$\mathbf{W}^{(l)} \in \mathbb{R}^{n_l \times n_{l-1}}$ and $\mathbf{b}^{(l)} \in \mathbb{R}^{n_l}$ are the weight matrix for linear mapping and the bias vector of layer $l$, respectively, where $n_0 = d_{in}$ is the input dimension. 
They are trainable parameters to be optimized through training. 
For a fully-connected layer, the number of trainable parameters is $(n_{l-1}+1) n_l$.

\begin{figure}[htp]
\centering
\includegraphics[width=1\linewidth]{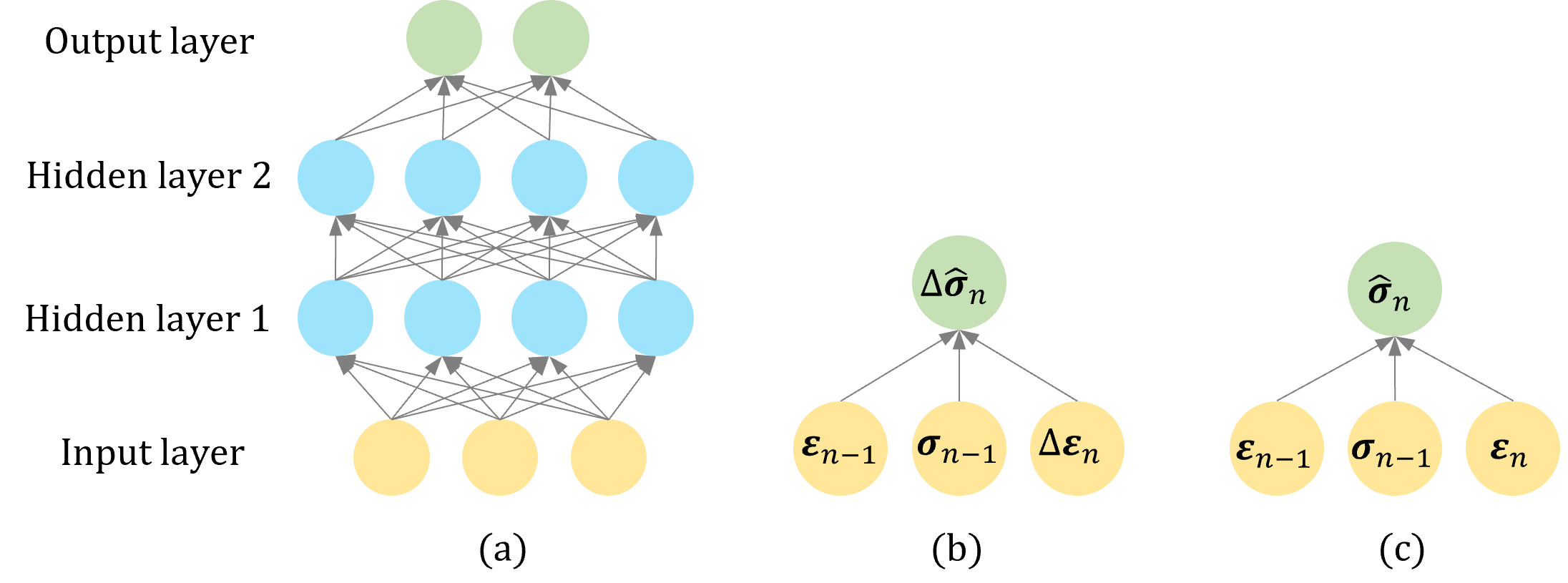}
\caption{Computational graphs of (a) a feed-froward deep neural network (DNN) with three input neurons, two hidden layers, and two output neurons, (b) a DNN constitutive model that takes one history step of stress-strain states and the current strain increment as input and predicts the current stress increment, and (c) a DNN constitutive model that takes one history step of stress-strain states and the current total strain as input and predicts the current total stress.}\label{fig.nn}
\end{figure}

Commonly used activation functions include the logistic sigmoid, the hyperbolic tangent function, the rectified linear unit (ReLU), and the leaky ReLU \cite{xu2015empirical}. 
Note that the choice of the activation of the output layer depends on the type of ML tasks. 
For regression tasks, which is the application of this study, a linear function is used in the output layer where the last hidden layer information is mapped to the output vector $\hat{\mathbf{y}}$, expressed as: $\hat{\mathbf{y}} = \mathbf{W}^{(L+1)} \mathbf{x}^{(L)} + \mathbf{b}^{(L+1)}$, where $\hat{\mathbf{y}}$ denotes the DNN approximation of the data output $\mathbf{y}$.
Fig. \ref{fig.nn}(a) shows the computational graph of a feed-forward DNN with three input neurons, two hidden layers, and two output neurons.

For path-dependent materials, the current stress response depends on the past stress-strain history. 
Fig. \ref{fig.nn}(b) shows the computational graph of a DNN constitutive model for path-dependent materials that takes one history step of stress-strain states and the current strain increment as input and predicts the current stress increment. 
Alternatively, the current strain increment can be replaced with the current total strain as input and the current total stress as the output, as shown in Fig. \ref{fig.nn}(c). 
The effects of the input/output representation will be discussed in the next subsection. 
These DNN constitutive models can be extended to consider pre-defined ISVs as input in addition to the measurable material states. 
An additional DNN can be used to model the evolution of the ISVs \cite{masi2021thermodynamics}. 
However, the dependency on pre-defined ISVs limits its applications, especially when only the measurable material states of path-dependent behaviors are available, e.g., the soil example to be demonstrated in Section \ref{sec:result_soil}.

Note that for DNN constitutive models, the number of history input steps is fixed once the model architecture is determined, which means the number of history input steps used for testing must be the same as that used in training. Furthermore, it is difficult for DNNs with recurrent connections from the output of step $n$ to the input of step $n+1$ to capture the essential information about the past history since the outputs are explicitly trained only to match the training set targets not being informed of the past history \cite{goodfellow2016deep}. These issues are addressed by RNNs introduced in the next subsection.

\subsection{Recurrent Neural Networks (RNNs) Constitutive Models}\label{sec:rnn-rnn}
\begin{figure}[htp]
\centering
\includegraphics[width=1\linewidth]{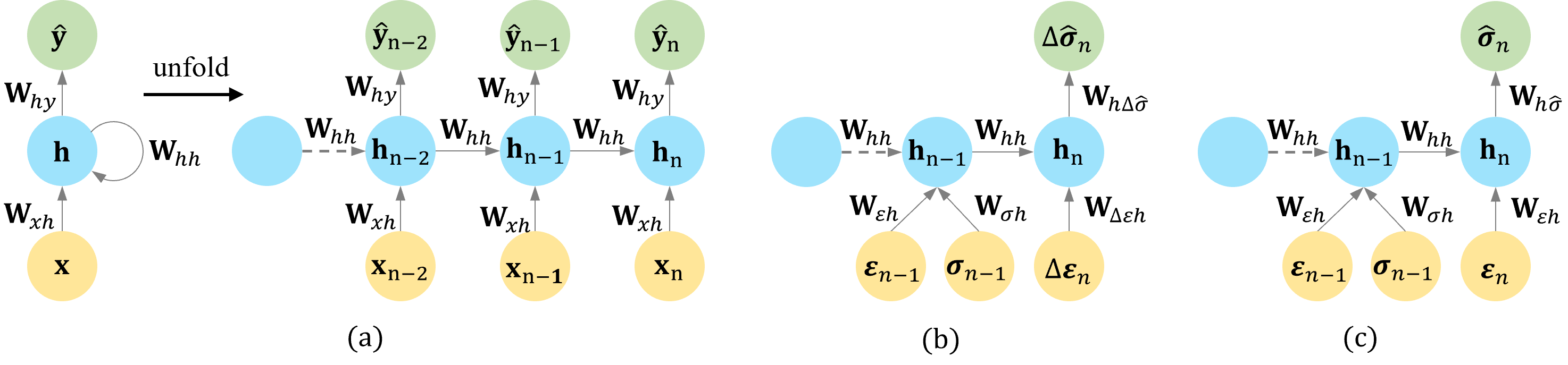}
\caption{Computational graphs of (a) a recurrent neural network (RNN), (b) a RNN constitutive model that takes one history step of stress-strain states and the current strain increment as input and predicts the current stress increment, and (c) an RNN constitutive model that takes one history step of stress-strain states and the current total strain as input and predicts the current total stress.}\label{fig.rnn}
\end{figure}

Recurrent neural networks (RNNs) designed for sequence learning have demonstrated successful applications in various domains, such as machine translation and speech recognition, due to their capability of learning history-dependent features that are essential for sequential prediction \cite{lipton2015critical,yu2019review}. 
Fig. \ref{fig.rnn}(a) illustrates the computational graphs of a folded RNN and an unfolded RNN, where $\mathbf{h}$ is a hidden state that captures essential history-dependent features from past information, which makes RNNs particularly suitable for modeling path-dependent material behaviors.
Unfolding of the RNN computational graph results in \textit{parameter sharing} across the network structure, reducing the number of trainable parameters and thus leading to more efficient training.
The length of input/output sequences can be arbitrary, which allows generalization to sequence lengths not appeared in the training set.
Each step can be viewed as a state. Despite the history sequence length, the trained RNN model always has the same input size, since it is specified in terms of transition from one state to another rather than in terms of a variable-length history of states \cite{goodfellow2016deep}. 
The forward propagation of RNN begins with an initial hidden state and the propagation equations at time step (state) $n$ are defined as
\begin{subequations}\label{eq.rnn-forward}
    \begin{align}
        \mathbf{h}_n & = a_{tanh} \big( \mathbf{W}_{hh} \mathbf{h}_{n-1} + \mathbf{W}_{xh} \mathbf{x}_{n} + \mathbf{b}_h \big), \\
        \hat{\mathbf{y}}_n & = \mathbf{W}_{hy} \mathbf{h}_n + \mathbf{b}_y,
    \end{align}
\end{subequations}
where $a_{tanh}$ is the hyperbolic tangent function; $\mathbf{W}_{xh}, \mathbf{W}_{hh}$, and $\mathbf{W}_{hy}$ are trainable weight coefficients for input-to-hidden, hidden-to-hidden, and hidden-to-output transformations, respectively; $\mathbf{b}_h$ and $\mathbf{b}_y$ are trainable bias coefficients. These trainable parameters are shared across all RNN steps.
Eqs. (\ref{eq.rnn-forward}a) transforms the previous hidden state $\mathbf{h}_{n-1}$ and the current input $\mathbf{x}_{n}$ to the current hidden state $\mathbf{h}_n$, while (\ref{eq.rnn-forward}b) transforms the current hidden state $\mathbf{h}_n$ to the current output $\hat{\mathbf{y}}_n$. The history information is captured by the hidden state of RNN by repeating the transformation in Eq. (\ref{eq.rnn-forward}a) for all RNN steps. The hidden state that carries the essential history-dependent information is passed to the final step and informs the final prediction.

Depending on applications, RNNs can have flexible architectures of input and output, such as one-to-one, one-to-many, many-to-one, and many-to-many \cite{goyal2018deep}. 
For example, the unfolded RNN shown in Fig. \ref{fig.rnn}(a) is a many-to-many type of RNN, which can be applied to, for example, name entity recognition. 
For path-dependent constitutive modeling, the many-to-one type of RNN is more suitable. 
Fig. \ref{fig.rnn}(b) illustrates the computational graph of an RNN constitutive model that takes one history step of stress-strain states and the current strain increment as input and predicts the current stress increment, defined as an \textit{incremental} RNN, whereas Fig. \ref{fig.rnn}(c) show the computational graph of an RNN constitutive model that takes one history step of stress-strain states and the current total strain as input and predicts the current total stress, defined as a \textit{total-form} RNN. 
Unlike the standard RNNs that has the same input size for all time steps (states), the history-step input of the RNN constitutive models shown in Fig. \ref{fig.rnn}(b)-(c) contains both strain and stress components whereas the current-step input contains only strain components. 
Considering one history step, the forward propagation of a typical \textit{total-form} RNN is expressed as follows
\begin{subequations}\label{eq.rnn-constitutive}
    \begin{align}
        & \mathbf{h}_{n-1} = a_{tanh} \big( \mathbf{W}_{hh} \mathbf{h}_{n-2} + \mathbf{W}_{\varepsilon h} \boldsymbol{\varepsilon}_{n-1} + \mathbf{W}_{\sigma h} \boldsymbol{\sigma}_{n-1} + \mathbf{b}_h \big), \\
        & \mathbf{h}_{n} = a_{tanh} \big( \mathbf{W}_{hh} \mathbf{h}_{n-1} + \mathbf{W}_{\varepsilon h} \boldsymbol{\varepsilon}_{n} + \mathbf{b}_h \big), \\
        & \hat{\boldsymbol{\sigma}}_{n} = \mathbf{W}_{h \hat{\sigma}} \mathbf{h}_{n} + \mathbf{b}_y,
    \end{align}
\end{subequations}
where $\mathbf{W}_{hh}$, $\mathbf{W}_{\varepsilon h}$, $\mathbf{W}_{\sigma h}$, and $\mathbf{W}_{h \hat{\sigma}}$ are trainable weight coefficients for hidden-to-hidden, strain-to-hidden, stress-to-hidden, and hidden-to-output transformations, respectively, and $\mathbf{b}_h$ and $\mathbf{b}_y$ are trainable bias coefficients.

To capture complex history-dependent patterns, deep RNNs are more advantageous \cite{goodfellow2016deep}. Similar to DNNs, stacking of fully-connected hidden layers can be added to affine input-to-hidden, hidden-to-hidden, and hidden-to-output transformations.

\begin{remark}
    Our studies show that the \textit{total-form RNN} is less sensitive to the strain increment than the \textit{incremental RNN}, as a consequence of \textbf{interpolation} outperforming \textbf{extrapolation}. 
    For instance, considering a training dataset with one stress-strain path and a constant loading (strain) increment, the final-step input (the strain increment) of the \textit{incremental RNN} is constant, whereas the final-step input (the total strain) of the \textit{total-form RNN} is not constant and covers the whole range of strain in the stress-strain path. 
    During testing on the same stress-strain path but with a different constant loading increment, larger or smaller than that of the training data, the \textit{incremental RNN} becomes inaccurate since the final-step input (the strain increment) of the testing data is beyond the range of the training strain increment and the prediction is an \textbf{extrapolation}. 
    In contrast, the \textit{total-form RNN} remains accurate because the final-step input (the total strain) of the testing data is within the range of training total strain and the prediction is an \textbf{interpolation}. Therefore, the proposed data-driven models in this study are built upon the total form.
\end{remark}

\subsubsection{Gated Recurrent Units (GRUs)}\label{sec:rnn-gru}
Standard RNNs suffer from short-term memory due to vanishing and exploding gradient issues that arise from recurrent connections \cite{bengio1993problem, bengio1994learning, goodfellow2016deep}. 
More effective RNNs for learning long-term dependencies have been developed, including the long short-term memory (LSTM) \cite{hochreiter1997long} cells and gated recurrent units (GRUs) \cite{cho2014properties, chung2014empirical}. 
A typical GRU consists of a \textit{reset} gate $\mathbf{r}_n$ that removes irrelevant past information, an \textit{update} gate $\mathbf{u}_n$ that controls the amount of past information passing to the next step, and a \textit{candidate} hidden state $\Tilde{\mathbf{h}}_n$ \cite{chung2014empirical}. Compared to LSTM, the GRU has fewer parameters as it lacks an \textit{output} gate. 
Considering one history step, the forward propagation of a typical GRU is expressed as follows
\begin{subequations}\label{eq.rnn-gru}
    \begin{align}
        \mathbf{r}_n & = a_{\sigma} \big( \mathbf{W}_{hr} \mathbf{h}_{n-1} + \mathbf{W}_{xr} \mathbf{x}_{n} + \mathbf{b}_r \big), \\
        \mathbf{u}_n & = a_{\sigma} \big( \mathbf{W}_{hu} \mathbf{h}_{n-1} + \mathbf{W}_{xu} \mathbf{x}_{n} + \mathbf{b}_u \big), \\
        \Tilde{\mathbf{h}}_n & = a_{tanh} \big( \mathbf{r}_n * \mathbf{W}_{h\Tilde{h}} \mathbf{h}_{n-1} + \mathbf{W}_{x\Tilde{h}} \mathbf{x}_{n} + \mathbf{b}_{\Tilde{h}} \big), \\
        \mathbf{h}_n & = \mathbf{u}_n * \mathbf{h}_{n-1} + (1-\mathbf{u}_n) * \Tilde{\mathbf{h}}_n + \mathbf{b}_{h}, \\
        \hat{\mathbf{y}}_n & = \mathbf{W}_{hy} \mathbf{h}_n + \mathbf{b}_{y},
    \end{align}
\end{subequations}
where $*$ denotes the Hadamard (element-wise) product; $a_{\sigma}$ is the sigmoid function; $a_{tanh}$ is the hyperbolic tangent function; $\mathbf{W}_{hr}$, $\mathbf{W}_{xr}$, $\mathbf{W}_{hu}$, $\mathbf{W}_{xu}$, $\mathbf{W}_{h\Tilde{h}}$, $\mathbf{W}_{x\Tilde{h}}$, and $\mathbf{W}_{hy}$ are trainable weight coefficients; $\mathbf{b}_r$, $\mathbf{b}_u$, $\mathbf{b}_{\Tilde{h}}$, $\mathbf{b}_{h}$, and $\mathbf{b}_{y}$ are trainable bias coefficients. Eq. (\ref{eq.rnn-gru}d) calculates the current hidden state $\mathbf{h}_n$ by a linear interpolation between the previous hidden state $\mathbf{h}_{n-1}$ and the candidate hidden state $\Tilde{\mathbf{h}}_n$, based on the update gate $\mathbf{u}_n$.
The RNN-based constitutive models proposed in this study are applicable to all types of RNNs for complicated path-dependent material behaviors with long-term dependent features. 

\subsubsection{Model Training}\label{sec:rnn-training}
Since the forward propagation (Eqs. \eqref{eq.rnn-forward}-\eqref{eq.rnn-gru}) is inherently sequential, i.e., each time step can only be computed after the previous one, the computation of the gradient of the loss function with respect to the trainable parameters cannot be parallelized and it needs to follow the reverse unfolded computational graph. The back-propagation through time algorithm is applied to RNNs \cite{goodfellow2016deep}.

\begin{figure}[htp]
\centering
\includegraphics[width=0.5\linewidth]{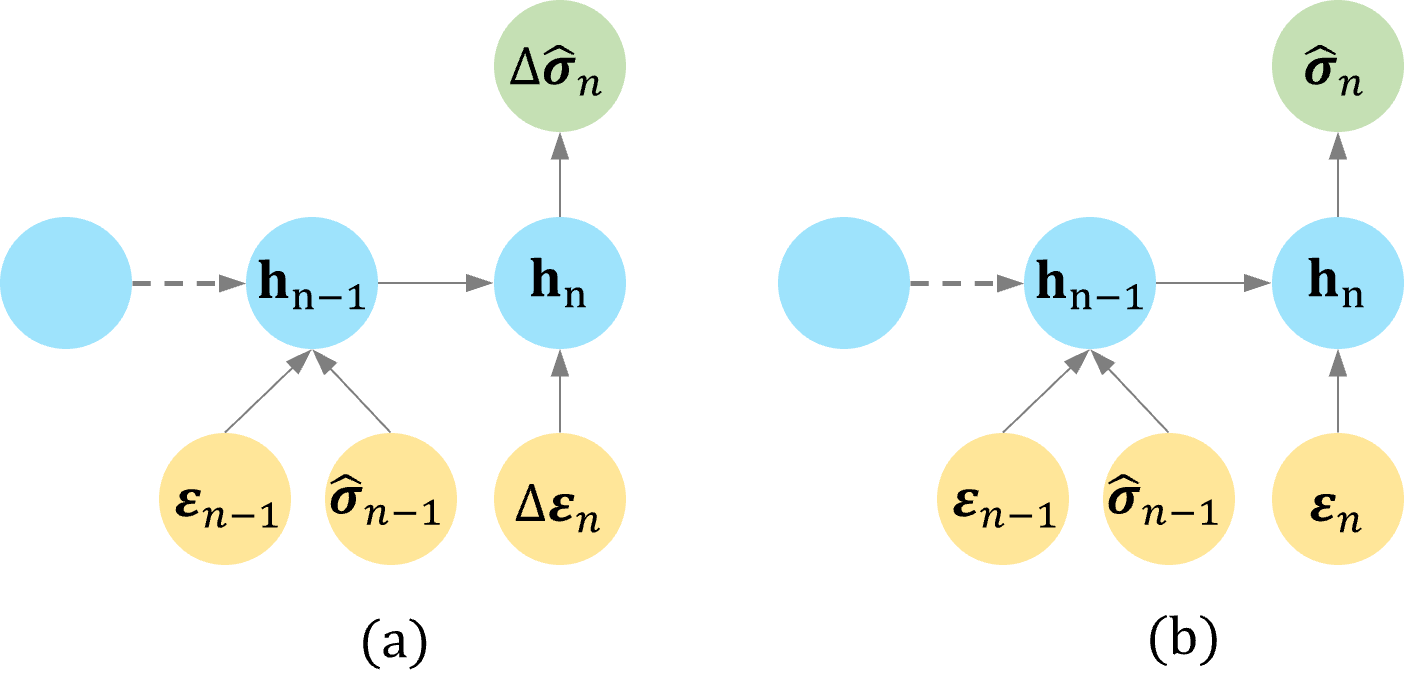}
\caption{Test mode of (a) a RNN constitutive model that takes one step history of stress-strain states and the current strain increment as input and predicts the current stress increment, and (b) an RNN constitutive model that takes one history step of stress-strain states and the current total strain as input and predicts the current total stress. The stress prediction from previous step is used as a part of the history input for the current-step prediction.}\label{fig.rnn_test}
\end{figure}

During training, the model receives the ground truth stress data of history steps, which is a \textit{teacher forcing} procedure emerging from the maximum likelihood criterion \cite{goodfellow2016deep}. 
However, the disadvantage of teacher forcing training arises when the trained model is applied in an \textit{open-loop} test mode with the network's previous outputs fed back as input for future predictions. 
The computational graphs of the test mode are shown in Fig. \ref{fig.rnn_test} corresponding to the RNN models (in the train mode) shown in Fig. \ref{fig.rnn}(b)-(c). 
In this case, the inputs the trained model receives could be quite different from those received during training, forcing the model to perform \textit{extrapolative} predictions and thus lead to large errors. 
Furthermore, such prediction errors could occur at the very first prediction, accumulate and propagate quickly, and contaminate the subsequent predictions. 
To mitigate the issue of error propagation due to the teacher forcing training and enhance model accuracy and robustness, we introduce \textit{stochasticity} to the training set by adding random perturbations to the ground truth stress data. In this way, the network can also learn the variability of the input conditions resembling those in the test mode.

\section{Thermodynamically Consistent Machine-Learned Internal State Variable Approach for Path-Dependent Materials}\label{sec:tcrnn}

\subsection{Thermodynamically Consistent Recurrent Neural Networks (TCRNNs)}
\begin{figure}[!h]
\centering
    \begin{subfigure}{0.495\textwidth}
        \centering
        \includegraphics[width=1\linewidth]{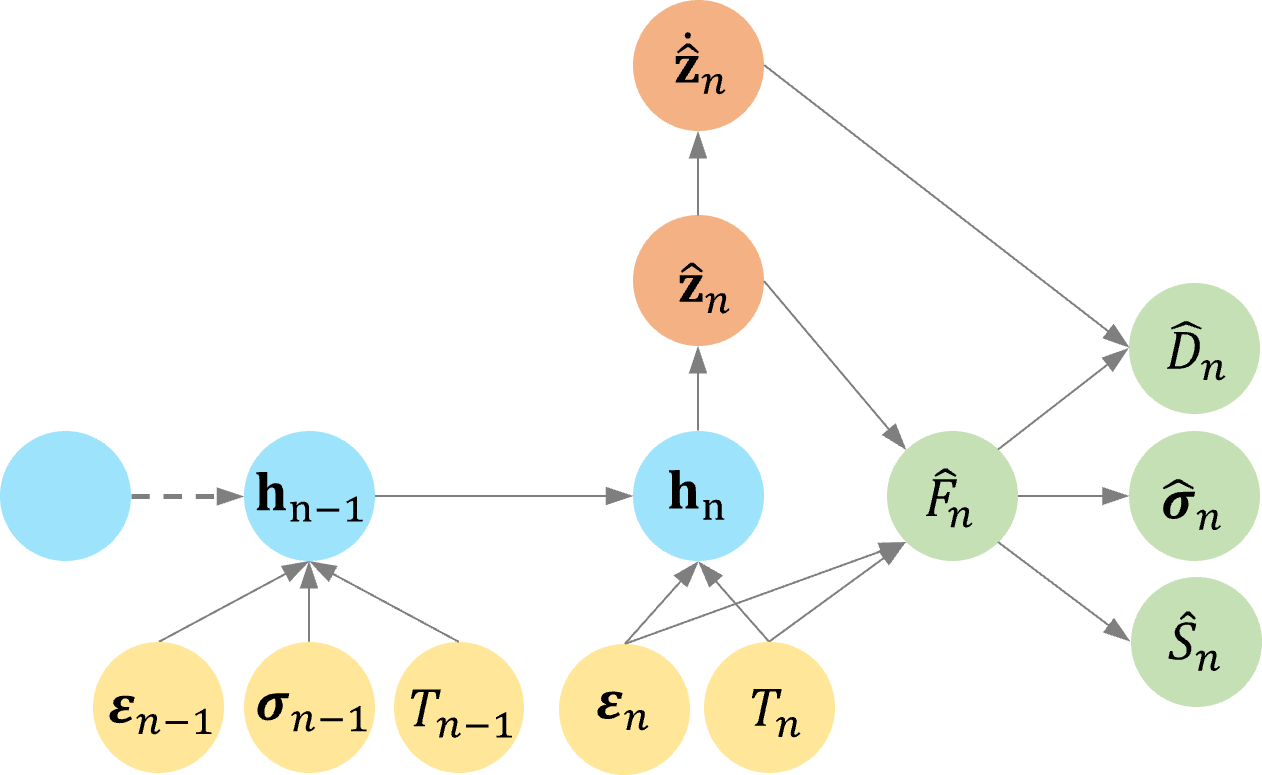}
        \caption{}
    \end{subfigure}
    \begin{subfigure}{0.495\textwidth}
        \centering
        \includegraphics[width=0.95\linewidth]{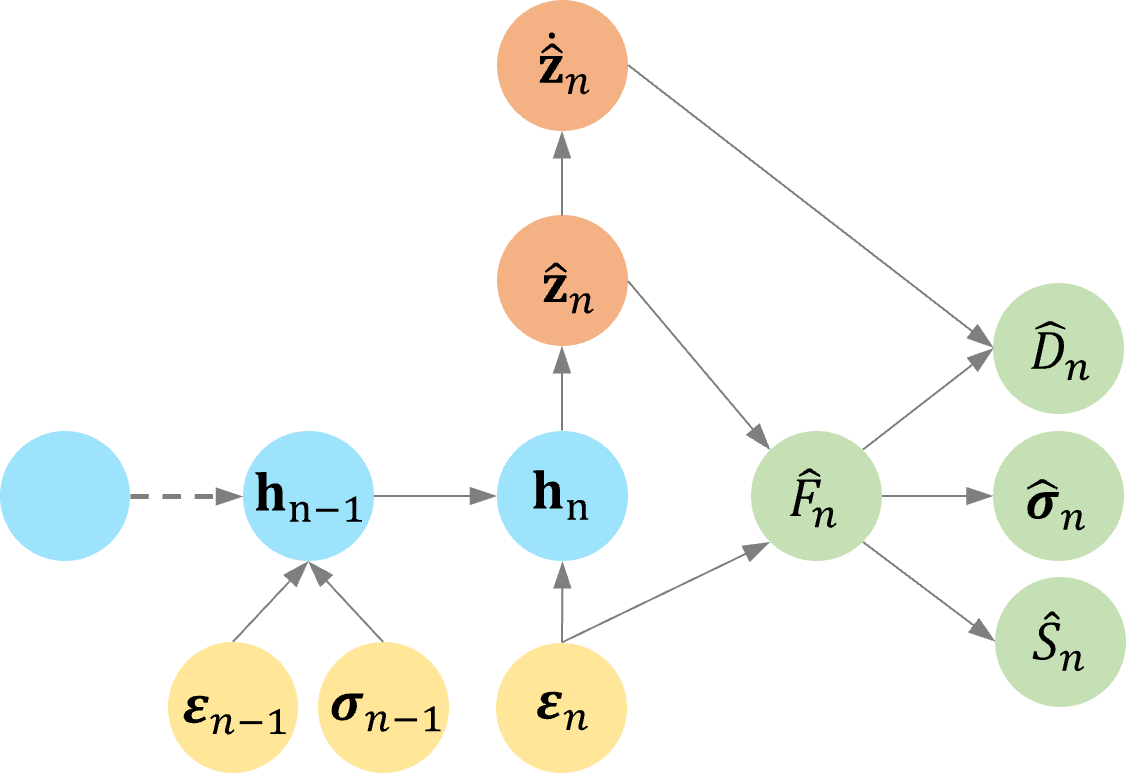}
        \caption{}
    \end{subfigure}
\caption{Computational graphs of (a) a thermodynamically consistent recurrent neural network (TCRNN) for non-isothermal processes and (b) a TCRNN for isothermal processes. The rate of the machine-learned internal state variable, $\dot{\hat{\mathbf{z}}}_n$, is computed and used for calculating the dissipation rate.}\label{fig.tcrnn_dzdt}
\end{figure}

To ensure thermodynamical consistency in the data-driven path-dependent materials modeling, thermodynamics principals introduced in Section \ref{sec:thermodynamics} are embedded into RNNs to extract the hidden ISVs. The proposed TCRNN consists of an RNN and a DNN. The computational graphs for non-isothermal or isothermal processes are illustrated in Fig. \ref{fig.tcrnn_dzdt}.

Since the hidden state $\mathbf{h}$ of RNNs captures essential history-dependent features from past material information, we propose to extract ISVs of materials from the hidden state of an RNN, as shown in Fig. \ref{fig.tcrnn_dzdt}. Considering one history step, the RNN-inferred ISV is expressed as follows
\begin{equation}\label{eq.tcrnn_rnn}
    \hat{\mathbf{z}}_n = f_{RNN}\big( \boldsymbol{\varepsilon}_{n}, T_{n}, \boldsymbol{\varepsilon}_{n-1}, \boldsymbol{\sigma}_{n-1}, T_{n-1} \big).
\end{equation}
Hereafter, a hat symbol is used to denote the predicted quantities. Considering one single layer for input-to-hidden, hidden-to-hidden, and hidden-to-output transformations in RNN, the RNN function $f_{RNN}$ consists of the following
\begin{subequations}\label{eq.rnn_eq}
    \begin{align}
        & \mathbf{h}_{n-1} = a_h\big( \mathbf{W}_{hh} \mathbf{h}_{n-2} + \mathbf{W}_{\varepsilon h} \boldsymbol{\varepsilon}_{n-1} + \mathbf{W}_{\sigma h} \boldsymbol{\sigma}_{n-1} + \mathbf{W}_{Th} T_{n-1} + \mathbf{b}_h \big), \\
        & \mathbf{h}_{n} = a_h\big( \mathbf{W}_{hh} \mathbf{h}_{n-1} + \mathbf{W}_{\varepsilon h} \boldsymbol{\varepsilon}_{n} + \mathbf{W}_{Th} T_{n} + \mathbf{b}_h \big), \\
        & \hat{\mathbf{z}}_{n} = a_z\big( \mathbf{W}_{hz} \mathbf{h}_{n} + \mathbf{b}_z \big),
    \end{align}
\end{subequations}
where $\mathbf{W}_{hh}$, $\mathbf{W}_{\varepsilon h}$, $\mathbf{W}_{\sigma h}$, $\mathbf{W}_{Th}$, and $\mathbf{W}_{hz}$ are trainable weight coefficients for hidden-to-hidden, strain-to-hidden, stress-to-hidden, temperature-to-hidden, and hidden-to-ISV transformations, respectively; $\mathbf{b}_h$ and $\mathbf{b}_z$ are trainable bias coefficients; $a_h(\cdot)$ and $a_z(\cdot)$ are activation functions. Note that the trainable parameters are shared across all steps of the RNN, which enhances training efficiency. $\hat{\mathbf{z}}_{n}$ is the machine-learned ISVs from the hidden state $\mathbf{h}_{n}$ of the RNN and its rate, $\dot{\hat{\mathbf{z}}}_{n}$, can be computed by automatic differentiation \cite{baydin2018automatic}, which will be discussed in the next subsection. 
Eqs. (\ref{eq.rnn_eq}a-b) transform the history and current measurable material states to the current hidden state $\mathbf{h}_{n}$ that carries the essential past information. For an RNN with more than one history step, Eq. (\ref{eq.rnn_eq}a) is repeated for all history steps. Eq. (\ref{eq.rnn_eq}c) infers the current ISV $\hat{\mathbf{z}}_{n}$ from the current hidden state $\mathbf{h}_{n}$.

\begin{figure}[!h]
    \centering
    \includegraphics[width=0.6\linewidth]{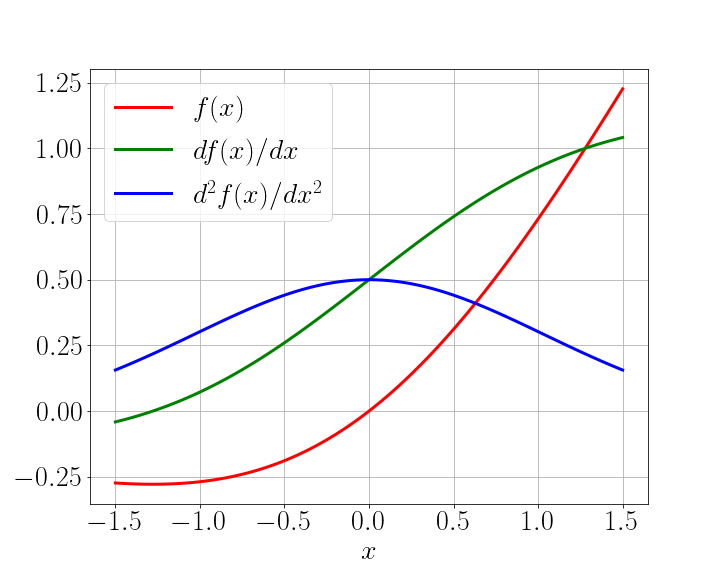}
    \caption{The SiLU activation function and its first-order and second-order gradients.}\label{fig.silu}
\end{figure}

A linear activation ($a_z(\cdot)$) is used for the transformation from the hidden state to the ISV in Eq. (\ref{eq.rnn_eq}c). The selection of the activation $a_h(\cdot)$ requires particular attention due to the issue of second-order vanishing gradients \cite{masi2021thermodynamics}. For effective training via back-propagation, the gradient of the output derivative with respect to the trainable parameters requires non-zero second-order gradients of activation functions. As a result, the activation function SiLU$(x)=x/(1+e^{-x})$ is selected for $a_h(\cdot)$ due to its smoothness and non-zero second-order gradients, as shown in Fig. \ref{fig.silu}.

Following the extraction of the ISV, a DNN is then used to predict the Helmholtz free energy given the strain, the temperature and the machine-learned ISVs, 
\begin{equation}\label{eq.tcrnn_nn}
    \hat{F}_n = f_{DNN}\big( \boldsymbol{\varepsilon}_{n}, T_{n}, \hat{\mathbf{z}}_{n} \big),
\end{equation}
where $f_{DNN}$ represents a DNN, as discussed in Section \ref{sec:rnn-nn}. The activation in the output layer is linear. The output Helmholtz free energy is then used to compute the stress $\hat{\boldsymbol{\sigma}}$, the dissipation rate $\hat{D}$, and the entropy $\hat{S}$ according to Eq. (\ref{eq.S_sigma_D}), which implicitly enforces the \textit{first thermodynamics principle}. The \textit{second thermodynamics principle}, i.e., $\hat{D} \ge 0$, is enforced in the loss function by constraining the network parameters, which will be discussed in the next subsection. The gradients of the output with respect to the inputs are obtained by automatic differentiation \cite{baydin2018automatic}. Since the output derivatives are involved in the loss function, SiLU is used for the activation of the hidden layers to avoid the issue of second-order vanishing gradients as discussed above. 

\subsection{Secondary Outputs}
To balance feature contributions to the loss function and accelerate the training process, the training dataset is standardized to have zero mean and unit variance. For instance, a feature $\mathbf{x}$ of the dataset is standardized by its mean $\mu_{\mathbf{x}}$  and standard deviation $std_{\mathbf{x}}$,
\begin{equation}\label{eq.standardization}
    \Bar{\mathbf{x}} = \frac{\mathbf{x} - \mu_{\mathbf{x}}}{std_{\mathbf{x}}}.
\end{equation}
Hereafter, a bar symbol is used to denote standardized quantities. From Eq. (\ref{eq.standardization}), we have
\begin{equation}\label{eq.standardization_grad}
    d\Bar{\mathbf{x}} = \frac{1}{std_{\mathbf{x}}} d\mathbf{x}.
\end{equation}

Considering standardized variables, Eqs. (\ref{eq.tcrnn_rnn}) and (\ref{eq.tcrnn_nn}) become
\begin{equation}\label{eq.tcrnn_rnn_bar}
    \hat{\mathbf{z}}_n = f_{RNN}\big( \bar{\boldsymbol{\varepsilon}}_{n}, \bar{T}_{n}, \bar{\boldsymbol{\varepsilon}}_{n-1}, \bar{\boldsymbol{\sigma}}_{n-1}, \bar{T}_{n-1} \big),
\end{equation}
\begin{equation}\label{eq.tcrnn_nn_bar}
    \hat{\bar{F}}_n = f_{DNN}\big( \bar{\boldsymbol{\varepsilon}}_{n}, \bar{T}_{n}, \hat{\mathbf{z}}_{n} \big).
\end{equation}
Therefore, the predicted stress is calculated by 
\begin{equation}\label{eq.stress_chain_rule}
    \hat{\boldsymbol{\sigma}}_n = \frac{\partial \hat{F}_n}{\partial \boldsymbol{\varepsilon}_n} = 
    \frac{\partial \hat{F}_n}{\partial \hat{\bar{F}}_n} 
    \frac{\partial \hat{\bar{F}}_n}{\partial \bar{\boldsymbol{\varepsilon}}_{n}} :
    \frac{\partial \bar{\boldsymbol{\varepsilon}}_{n}}{\partial \boldsymbol{\varepsilon}_n} = 
    \frac{std_F}{std_{\varepsilon}} \frac{\partial \hat{\bar{F}}_n}{\partial \bar{\boldsymbol{\varepsilon}}_{n}}.
\end{equation}
with 
\begin{equation}\label{eq.std}
    \frac{\partial \hat{F}_n}{\partial \hat{\bar{F}}_n} = std_F 
    \quad \text{and} \quad
    \frac{\partial \bar{\boldsymbol{\varepsilon}}_{n}}{\partial \boldsymbol{\varepsilon}_n} = \frac{1}{std_{\varepsilon}} \mathbf{I},
\end{equation}
according to Eq. (\ref{eq.standardization_grad}). $\mathbf{I}$ is the second-order identity tensor. $\frac{\partial \hat{\bar{F}}_n}{\partial \bar{\boldsymbol{\varepsilon}}_{n}}$ is the gradient of the output with respect to the input in Eq. (\ref{eq.tcrnn_nn_bar}) and can be obtained by automatic differentiation \cite{baydin2018automatic}.
Similarly, the predicted entropy is computed by
\begin{equation}\label{eq.entropy_chain_rule}
    \hat{S}_n = - \frac{\partial \hat{F}_n}{\partial T_n} = 
    - \frac{\partial \hat{F}_n}{\partial \hat{\bar{F}}_n} 
    \frac{\partial \hat{\bar{F}}_n}{\partial \bar{T}_{n}} 
    \frac{\partial \bar{T}}{\partial T_n} = 
    - \frac{std_F}{std_T} \frac{\partial \hat{\bar{F}}_n}{\partial \bar{T}_{n}}.
\end{equation}
The predicted dissipation rate is computed by
\begin{equation}\label{eq.dissipation_chain_rule}
    \hat{D}_n = - \frac{\partial \hat{F}_n}{\partial \hat{\mathbf{z}}_n} \cdot \dot{\hat{\mathbf{z}}}_n = 
    - \frac{\partial \hat{F}_n}{\partial \hat{\bar{F}}_n} 
    \frac{\partial \hat{\bar{F}}_n}{\partial \hat{\mathbf{z}}_{n}} \cdot \dot{\hat{\mathbf{z}}}_n = 
    - std_F \frac{\partial \hat{\bar{F}}_n}{\partial \hat{\mathbf{z}}_{n}} 
    \cdot \dot{\hat{\mathbf{z}}}_n,
\end{equation}
where $\dot{\hat{\mathbf{z}}}_n$ can be obtained by applying the chain rule to Eq. (\ref{eq.tcrnn_rnn_bar})
\begin{subequations}\label{eq.zdot_chain_rule}
    \begin{align}
        \dot{\hat{\mathbf{z}}}_n & = 
        \frac{\partial \hat{\mathbf{z}}_n}{\partial \bar{\boldsymbol{\varepsilon}}_n}: \dot{\bar{\boldsymbol{\varepsilon}}}_n + 
        \frac{\partial \hat{\mathbf{z}}_n}{\partial \bar{T}_n} \dot{\bar{T}}_n +
        \frac{\partial \hat{\mathbf{z}}_n}{\partial \bar{\boldsymbol{\varepsilon}}_{n-1}}: \dot{\bar{\boldsymbol{\varepsilon}}}_{n-1} +
        \frac{\partial \hat{\mathbf{z}}_n}{\partial \bar{\boldsymbol{\sigma}}_{n-1}}: \dot{\bar{\boldsymbol{\sigma}}}_{n-1} +
        \frac{\partial \hat{\mathbf{z}}_{n-1}}{\partial \bar{T}_n} \dot{\bar{T}}_{n-1}, \\
        & = 
        \frac{\partial \hat{\mathbf{z}}_n}{\partial \bar{\boldsymbol{\varepsilon}}_n}: \frac{\dot{\boldsymbol{\varepsilon}}_n}{std_{\varepsilon}} + 
        \frac{\partial \hat{\mathbf{z}}_n}{\partial \bar{T}_n} \frac{\dot{T}_n}{std_T} +
        \frac{\partial \hat{\mathbf{z}}_n}{\partial \bar{\boldsymbol{\varepsilon}}_{n-1}}: \frac{\dot{\boldsymbol{\varepsilon}}_{n-1}}{std_{\varepsilon}} +
        \frac{\partial \hat{\mathbf{z}}_n}{\partial \bar{\boldsymbol{\sigma}}_{n-1}}: \frac{\dot{\boldsymbol{\sigma}}_{n-1}}{std_{\sigma}} +
        \frac{\partial \hat{\mathbf{z}}_n}{\partial \bar{T}_{n-1}} \frac{\dot{T}_{n-1}}{std_T},
    \end{align}
\end{subequations}
which requires the rate of the input variables, including $\dot{\boldsymbol{\varepsilon}}_n$, $\dot{T}_n$, $\dot{\boldsymbol{\varepsilon}}_{n-1}$, $\dot{\boldsymbol{\sigma}}_{n-1}$, $\dot{T}_{n-1}$, etc.

\begin{figure}[!h]
\centering
    \begin{subfigure}{0.495\textwidth}
        \centering
        \includegraphics[width=1\linewidth]{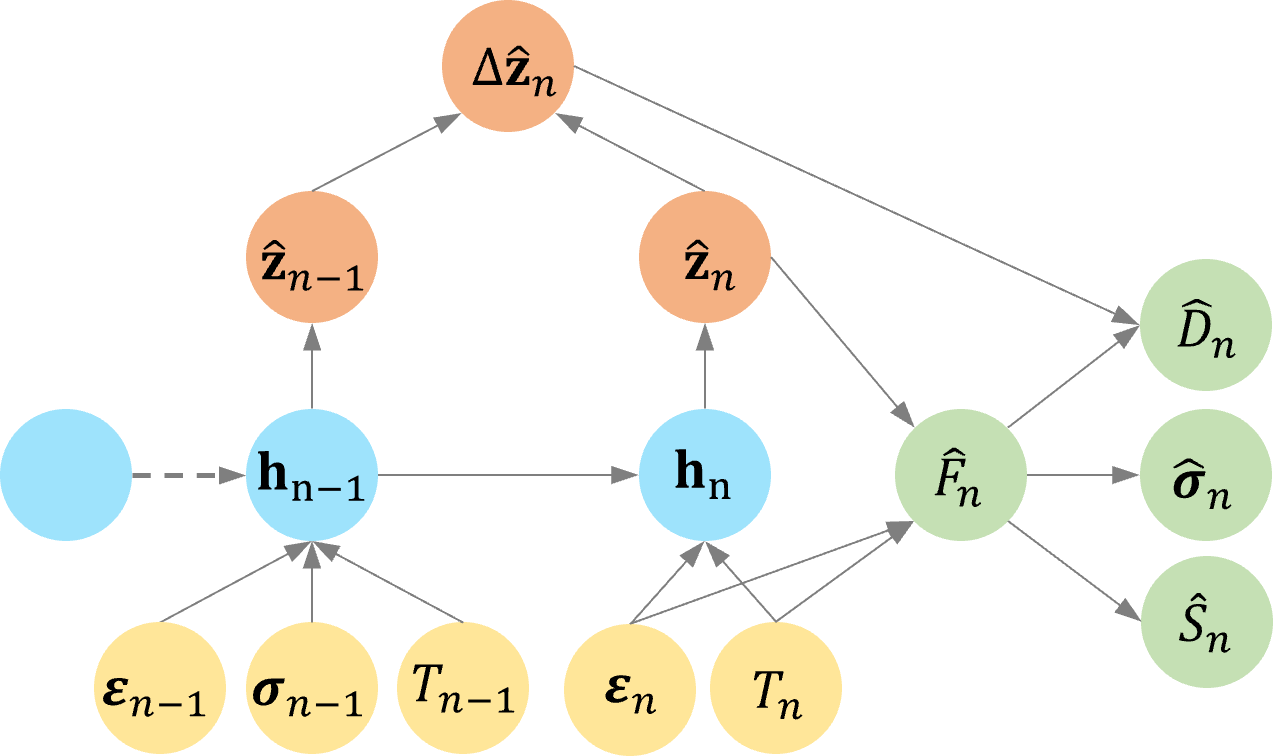}
        \caption{}
    \end{subfigure}
    \begin{subfigure}{0.495\textwidth}
        \centering
        \includegraphics[width=0.95\linewidth]{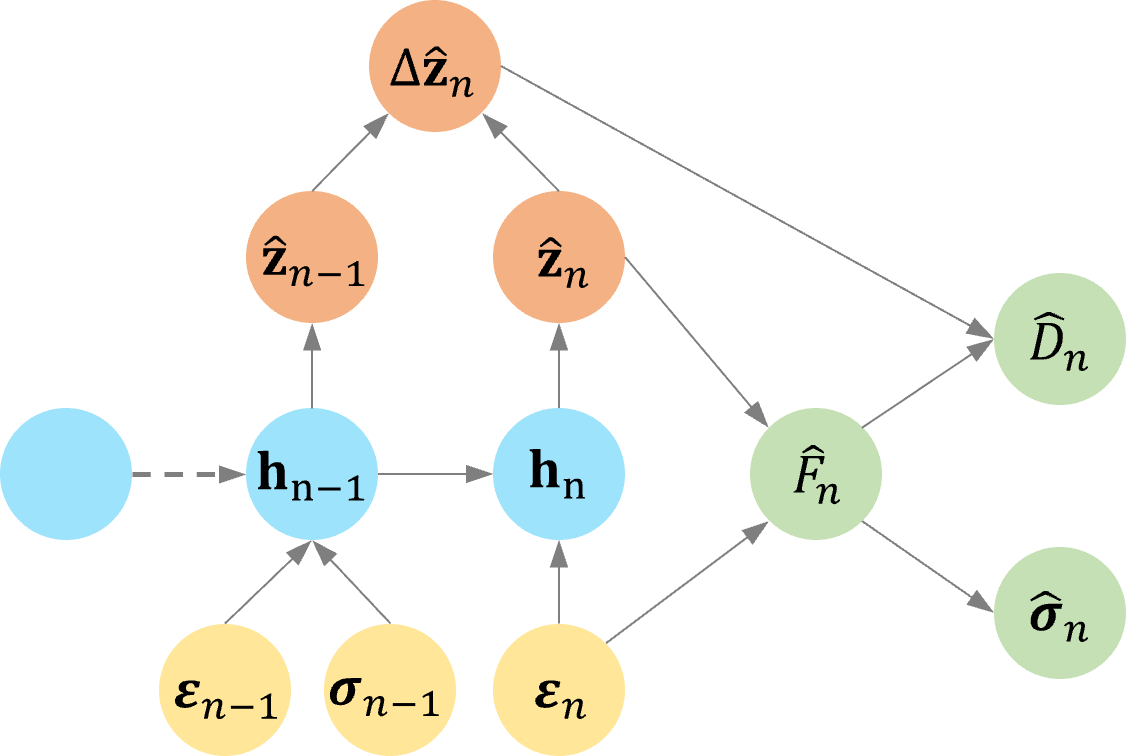}
        \caption{}
    \end{subfigure}
\caption{Computational graphs of (a) a thermodynamically consistent recurrent neural network (TCRNN) for non-isothermal processes and (b) a TCRNN for isothermal processes. The increment of the machine-extracted internal state, $\Delta \hat{\mathbf{z}}_n$, is computed and used for calculating the dissipation rate.}\label{fig.tcrnn_dz}
\end{figure}

The direct calculation of $\dot{\hat{\mathbf{z}}}_n$ through Eq. (\ref{eq.zdot_chain_rule}) can be computationally intractable, especially when the dataset size, the internal state dimension, and the number RNN input steps are large. 
Alternatively, the rate of the ISVs can be approximated by $\dot{\hat{\mathbf{z}}}_n \approx \Delta \hat{\mathbf{z}}_n/\Delta t$, where $\Delta \hat{\mathbf{z}}_n = \hat{\mathbf{z}}_n - \hat{\mathbf{z}}_{n-1}$ is the increment of the ISVs at the current step $n$. 
To obtain $\Delta \mathbf{z}_n$, alternative TCRNNs are proposed, as shown in Fig. \ref{fig.tcrnn_dz}, where the last two steps of the RNN infer the ISVs, $\hat{\mathbf{z}}_{n-1}$ and $\hat{\mathbf{z}}_n$, respectively. These TCRNN models enhance training efficiency.

\subsection{Model Training}
The loss function is expressed as
\begin{equation}\label{eq.tcrnn_loss1}
    Loss = \sum_n ||\bar{\boldsymbol{\sigma}}_n - \hat{\bar{\boldsymbol{\sigma}}}_n ||_{L_1}^2 + \beta_1 ||\bar{F}_n - \hat{\bar{F}}_n ||_{L_1}^2 + \beta_2 ||\bar{D}_n - \hat{\bar{D}}_n ||_{L_1}^2 + \beta_3 ||\bar{S}_n - \hat{\bar{S}}_n ||_{L_1}^2,
\end{equation}
where $\beta_i$, $i=1,..,3$, are regularization parameters. $||\cdot||_{L_1}$ denotes the $L_1$ norm, and $\beta_3$ is set to be zero if the data of the entropy is unavailable. The training consists of forward propagation and backward propagation. During the forward propagation, the machine-learned ISVs are implicitly embedded in the calculation of the predicted measurable quantities by following the thermodynamics principles. During the backward propagation, the errors of the measurable quantities are back-propagated to update the model's trainable parameters and refine machine-learned ISVs.

In some cases where the data of the dissipation rate $D$ is unavailable, the non-negativity condition, i.e., $\hat{D} \ge 0$, can be imposed instead, which is resulted from the thermodynamics second law, Eq. (\ref{eq.mech_diss_rate}). To this end, the rectified linear unit (ReLU) can be utilized, i.e., ReLU(x) = $max(0,x) \ge 0$, which is positive only if $x$ is positive. Hence, ReLU$(-\hat{D})$ is positive only if $\hat{D}$ is negative, which corresponds to violation of the non-negativity condition $\hat{D} \ge 0$. Including ReLU$(-\hat{D})$ to the loss function penalizes the violation of the non-negativity condition and enforces $\hat{D} \ge 0$ to be satisfied, which imposes a constraint on the network parameters during training. 
The loss function becomes
\begin{equation}\label{eq.tcrnn_loss2}
    Loss = \sum_n ||\bar{\boldsymbol{\sigma}}_n - \hat{\bar{\boldsymbol{\sigma}}}_n ||_{L_1}^2 + \beta_1 ||\bar{F}_n - \hat{\bar{F}}_n ||_{L_1}^2 + \beta_2 ReLU(- \hat{D}_n).
\end{equation}

Similarly, if the data of the Helmholtz free energy $F$ is unavailable, the non-negativity condition, i.e., $\hat{F} \ge 0$, can be imposed by adding ReLU$(-\hat{F}_n)$ to the loss function,
\begin{equation}\label{eq.tcrnn_loss3}
    Loss = \sum_n ||\bar{\boldsymbol{\sigma}}_n - \hat{\bar{\boldsymbol{\sigma}}}_n ||_{L_1}^2 + \beta_1 ReLU(-\hat{F}_n) + \beta_2 ReLU(- \hat{D}_n).
\end{equation}

In some cases where prior knowledge of certain ISVs are available, the TCRNN models can be trained in a \textit{hybrid} mode by leveraging the existing ISVs and simultaneously inferring additional thermodynamically consistent ISVs that are essential to path-dependent behaviors. Considering $\mathbf{z}^p_n$ as the known ISVs and $\bar{\mathbf{z}}^p_n$ as the corresponding standardized quantity, the loss function becomes
\begin{equation}\label{eq.tcrnn_loss4}
    Loss = \sum_n ||\bar{\boldsymbol{\sigma}}_n - \hat{\bar{\boldsymbol{\sigma}}}_n ||_{L_1}^2 + \beta_1 ReLU(-\hat{F}_n) + \beta_2 ReLU(- \hat{D}_n) + \beta_4||\bar{\mathbf{z}}^p_n - \hat{\bar{\mathbf{z}}}^p_n ||_{L_1}^2,
\end{equation}
where the last term enables the TCRNN model to learn the existing ISVs. For the TCRNN model to infer additional essential ISVs, the prescribed internal state dimension, $|\hat{\mathbf{z}}|$, should be greater than the dimension of the existing ISVs, $|\mathbf{z}^p|$. Note that both existing and machine-learned ISVs are passed to the DNN to predict the Helmholtz free energy (Eq. \ref{eq.tcrnn_nn_bar}) and the downstream calculations (Eqs. \ref{eq.dissipation_chain_rule}-\ref{eq.zdot_chain_rule}).

Apart from thermodynamics, the time (or self) consistency condition is critical for convergence of numerical approximation when $\Delta t \rightarrow 0$ \cite{xu2021learning}.
\begin{equation}\label{eq.time_consistency}
    \lim_{\Delta \boldsymbol{\varepsilon} \rightarrow 0} \Delta \hat{\boldsymbol{\sigma}} = 0.
\end{equation}
To achieve the time consistency condition, the training set can be augmented by additional samples constituted by zero strain increment and zero stress increment at different material states (time steps), which enables the machine-learned material model to learn the time consistency condition from data. Alternatively, the self-consistency condition can be integrated into the RNN architecture by definition \cite{bonatti2022importance}.

The optimal parameters of TCRNN are obtained by minimizing the loss functions (Eqs. (\ref{eq.tcrnn_loss1})-(\ref{eq.tcrnn_loss4})) using the open-source Pytorch library \cite{paszke2019pytorch}.
The Adam optimizer \cite{kingma2014adam} is adopted for back-propagation training and the initial learning rate is set to be 0.001. 
Since the training dataset is standardized to balance feature contributions to the loss function, the training is less sensitive to the regularization parameters. Hence, the regularization parameters are set to be 1 unless further tuning is required to achieve better training performances.
\section{Numerical Results}\label{sec:result}

\subsection{Modeling Elasto-Plastic Materials}\label{sec:result_ep}
To demonstrate the accuracy, robustness, and generalization performances, the proposed TCRNN is applied to model a material with synthetic data generated by the one-dimensional elasto-plastic material with kinematic hardening. The Helmholtz free energy potential is expressed as
\begin{equation}\label{eq.1Dep_F}
    F(\varepsilon,\varepsilon^p) = \frac{E}{2} (\varepsilon - \varepsilon^p)^2 + \frac{H}{2} (\varepsilon^p)^2,
\end{equation}
where $E=100$ $GPa$ is the Young's modulus; $H=100$ $GPa$ is the kinematic hardening parameter; $\varepsilon$ is the total strain; $\varepsilon^p$ is the plastic strain, which can be considered as a phenomenological ISV. The yield stress is $k=100$ $MPa$. The stress and the dissipation rate can be obtained by Eq. (\ref{eq.S_sigma_D}b-c) as follows
\begin{subequations}\label{eq.1Dep_s_D}
    \begin{align}
    \sigma & = \frac{\partial F}{\partial \varepsilon} = E(\varepsilon-\varepsilon^p), \\
    D & = - \frac{\partial F}{\partial \varepsilon^p} \dot{\varepsilon}^p = (\sigma-H\varepsilon^p) \dot{\varepsilon}^p.
    \end{align}
\end{subequations}

\begin{figure}[htp]
    \centering
    \includegraphics[width=0.4\linewidth]{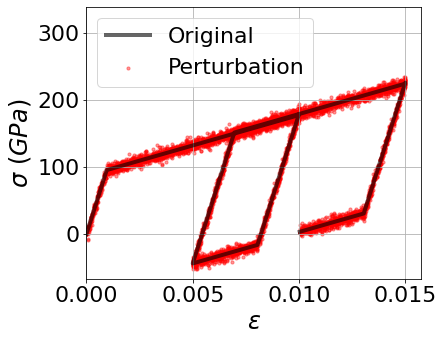}
    \caption{Training data with random perturbations to enhance prediction accuracy and robustness of TCRNNs, where the black solid line denotes the original data and red circles denote randomly perturbed data. $r=0.3$ is employed to generate the random perturbations.}\label{fig.1Dep_noisy_data}
\end{figure}

The dataset is generated by Eqs. (\ref{eq.1Dep_F})-(\ref{eq.1Dep_s_D}), which contains five samples with the same stress-strain path (two loading-unloading cycles), as show in Fig. \ref{fig.1Dep_noisy_data}. 
The only difference in these samples is the strain increment, ranging from $3.75 \times 10^{-5}$ to $7.5 \times 10^{-5}$. The data of the sample with a strain increment of $5.0 \times 10^{-5}$ is used to train the TCRNN. 
The remaining samples are in the testing set to evaluate the trained model.

To address the issue of error propagation in the test mode due to the teacher forcing training, as discussed in Section \ref{sec:rnn-training}, and enhance model accuracy and robustness, \textit{stochasticity} is introduced to the stress data so that the model learns the  uncertainties of the input conditions resembling those in the test mode. 
Random perturbations are generated from a normal (Gaussian) distribution with a zero mean and a standard deviation of $r \times \sigma_{max}$, where $\sigma_{max}$ is the maximum stress in the data and $r$ is a user-defined parameter to control the level of randomness. 
Fig. \ref{fig.1Dep_noisy_data} shows the original stress-strain data in a black solid line and the randomly perturbed data in red circles. During supervised training, the noiseless stress is the ground truth and the input stress variable is no longer deterministic but rather stochastic, contributed by the random perturbations. 

The TCRNN model \XH{based on the time rates of ISVs ($\dot{\hat{\mathbf{z}}}$)}, as shown in Fig. \ref{fig.tcrnn_dzdt}(b), is employed in this example. 
A GRU is used to infer the ISV $\hat{\mathbf{z}}$ and describe its evolution by following the thermodynamics second law.
\XH{The GRU consists of one hidden layer for all affine transformations in Eq. \eqref{eq.rnn-gru} with the model complexity represented by the dimension of the hidden state, $|\mathbf{h}|$.}
Since the data of the free energy $F$ and the dissipation rate $D$ can be obtained by Eq. (\ref{eq.1Dep_s_D}) in this example, the loss function in Eq. \eqref{eq.tcrnn_loss1} is employed with $\beta_1=\beta_2=1$ and $\beta_3=0$. 
A relative error used to measure the prediction accuracy is defined as follows
\begin{equation}\label{eq.relative_error}
    e = \frac{||\boldsymbol{\Sigma} - \hat{\boldsymbol{\Sigma}}||_{L_2}}{||\boldsymbol{\Sigma}||_{L_2}},
\end{equation}
where $\boldsymbol{\Sigma}$ and $\hat{\boldsymbol{\Sigma}}$ contain the stress data and stress predictions at all time steps of a loading path, respectively.

\XH{In the following subsections, the effects of the strain increments on model performance are first investigated. Since the machine-inferred ISVs are critical to the accuracy of path-dependent materials modeling, various factors that can affect the quality of the machine-inferred ISVs are investigated, including the number of RNN steps, the internal state dimension \XH{($|\hat{\mathbf{z}}|$)}, and the model complexity \XH{($|\mathbf{h}|$)}. Further, the generalization capabilities of the TCRNN model are examined.}

\subsubsection{Effects of Strain Increments}
\begin{figure}[htp]
\centering
    \begin{subfigure}{1\textwidth}
        \centering
        \includegraphics[width=1\linewidth]{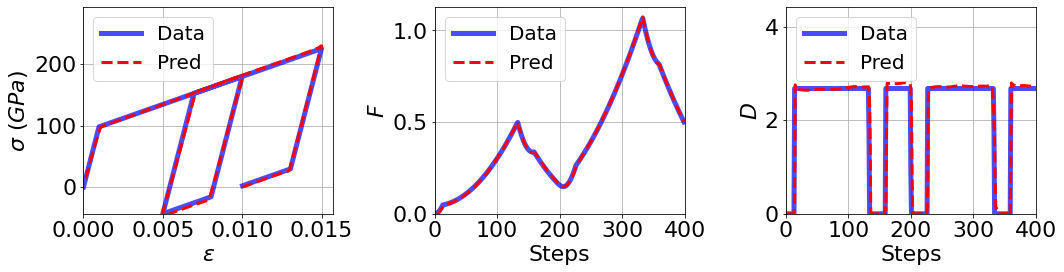}
        \vspace{-1.5\baselineskip}
        \caption{Testing: $|\Delta \varepsilon|=7.5 \times 10^{-5}$}
    \vspace{0.1in}
    \end{subfigure}
    \begin{subfigure}{1\textwidth}
        \centering
        \includegraphics[width=1\linewidth]{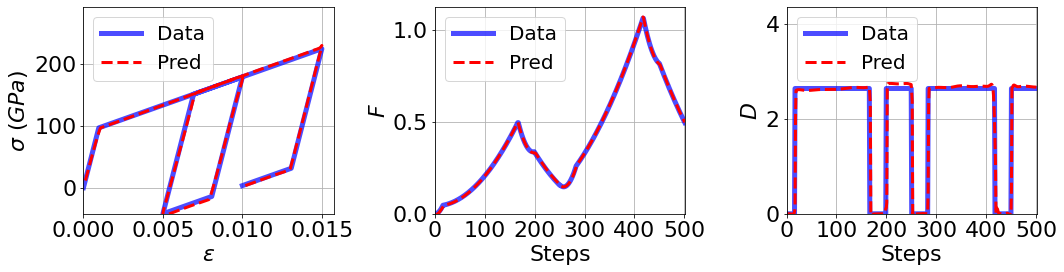}
        \vspace{-1.5\baselineskip}
        \caption{Testing: $|\Delta \varepsilon|=6.0 \times 10^{-5}$}
    \vspace{0.1in}
    \end{subfigure}
    \begin{subfigure}{1\textwidth}
        \centering
        \includegraphics[width=1\linewidth]{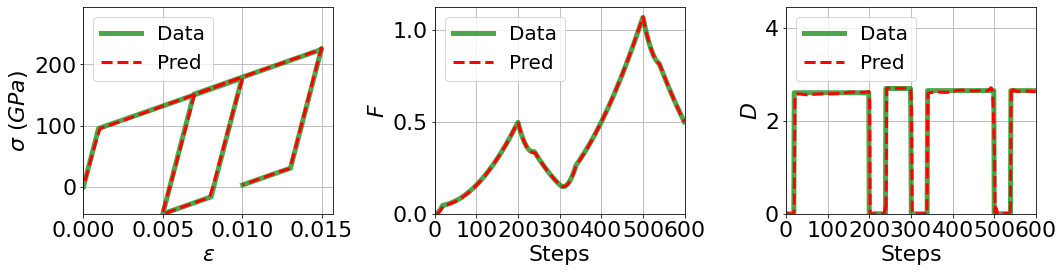}
        \vspace{-1.5\baselineskip}
        \caption{Training: $|\Delta \varepsilon|=5.0 \times 10^{-5}$}
    \vspace{0.1in}
    \end{subfigure}
    \begin{subfigure}{1\textwidth}
        \centering
        \includegraphics[width=1\linewidth]{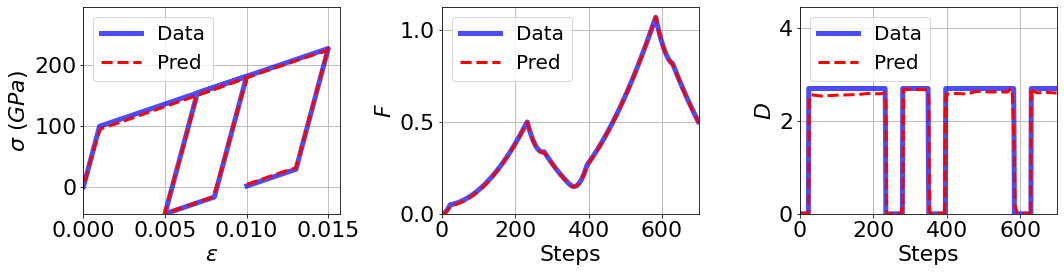}
        \vspace{-1.5\baselineskip}
        \caption{Testing: $|\Delta \varepsilon|=4.29 \times 10^{-5}$}
    \vspace{0.1in}
    \end{subfigure}
    \begin{subfigure}{1\textwidth}
        \centering
        \includegraphics[width=1\linewidth]{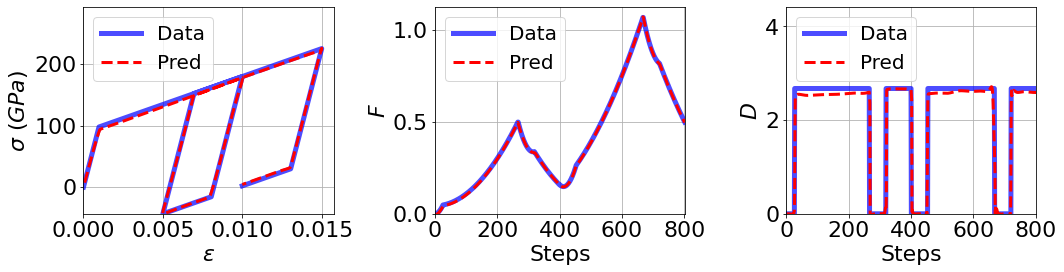}
        \vspace{-1.5\baselineskip}
        \caption{Testing: $|\Delta \varepsilon|=3.75 \times 10^{-5}$}
    \end{subfigure}
\caption{Comparison of predictions of the TCRNN with data: (a) testing case with a strain increment of $|\Delta \varepsilon|=7.5 \times 10^{-5}$; (b) testing case with a strain size of $|\Delta \varepsilon|=6.0 \times 10^{-5}$; (c) training case with a strain size of $|\Delta \varepsilon|=7.5 \times 10^{-5}$; (d) testing case with a strain size of $|\Delta \varepsilon|=4.29 \times 10^{-5}$; (e) testing case with a strain size of $|\Delta \varepsilon|=3.75 \times 10^{-5}$.}\label{fig.1Dep_load_increment}
\end{figure}

\begin{figure}[htp]
    \centering
    \includegraphics[width=0.4\linewidth]{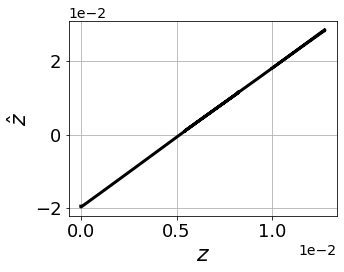}
    \caption{Correlation between the machine-learned internal state variable and the phenomenological internal state variable.}\label{fig.1Dep_z_vs_zhat}
\end{figure}

We first investigate how the strain increment affects the prediction accuracy of the TCRNN model. The model has a scalar ISV and \XH{a hidden state dimension of 30 ($|\mathbf{h}|=30$)}.
Fig. \ref{fig.1Dep_load_increment} compares the predictions with data, where the case with a green color line is used for training and those with blue color lines are used for testing. 
The trained model achieves 1.1$\%$  relative error (Eq. \eqref{eq.relative_error}) on the stress prediction for the training case and 1.9$\%$ mean relative error for the testing cases. 
The mean relative error is obtained by averaging the relative errors (Eq. \eqref{eq.relative_error}) of all cases.
The good agreement between predictions and data for all quantities, including the stress, the Helmholtz free energy, and the dissipation rate demonstrates that the model maintains high prediction accuracy and robustness as the strain increment varies.
Fig. \ref{fig.1Dep_z_vs_zhat} shows that the machine-learned ISV is monotonically correlated with the phenomenological ISV, which demonstrates the capability of the TCRNN model in extracting mechanistically and thermodynamically consistent ISV essential to dissipative elasto-plastic material behaviors.

\subsubsection{Effects of The Number of RNN Steps}
\begin{figure}[htp]
\centering
    \begin{subfigure}{0.49\textwidth}
        \centering
        \includegraphics[width=1\linewidth]{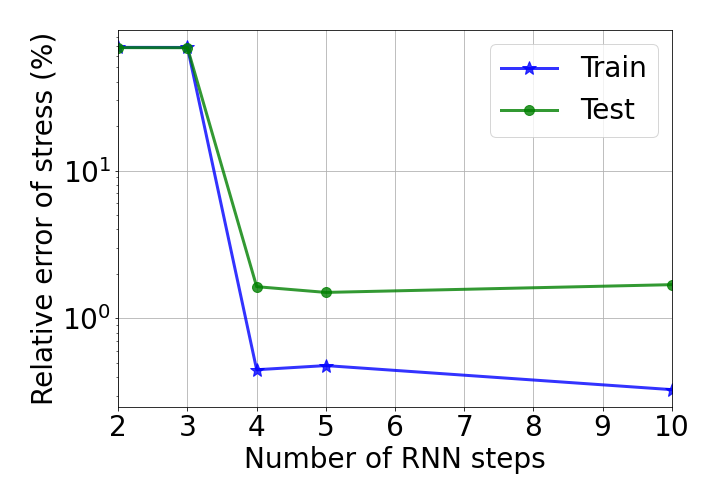}
        \vspace{-2\baselineskip}
        \caption{}
    \end{subfigure}
    \begin{subfigure}{0.49\textwidth}
        \centering
        \includegraphics[width=1\linewidth]{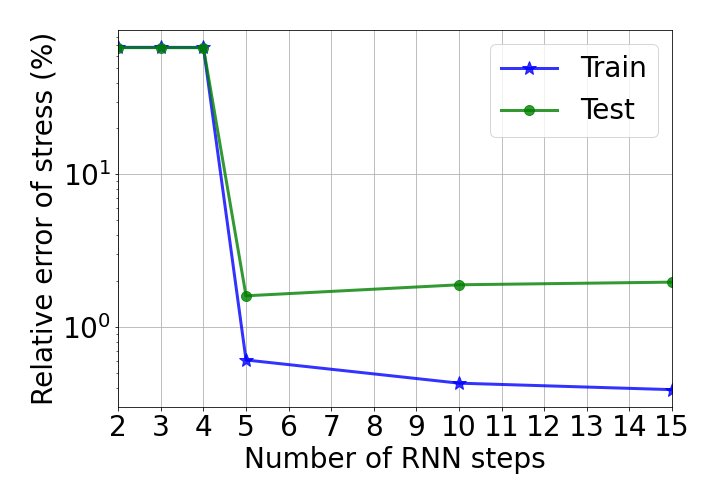}
        \vspace{-2\baselineskip}
        \caption{}
    \end{subfigure}
\caption{Effects of the number of RNN steps on the accuracy of models with a scalar internal state variable and: (a) \XH{a hidden state dimension $|\mathbf{h}|=20$}; (b) \XH{a hidden state dimension $|\mathbf{h}|=30$}.}\label{fig.1Dep_nsteps}
\end{figure}

In the second example, TCRNN models with a scalar ISV and various RNN steps (including history and current steps) are examined. \XH{Note that the model complexity is independent of the number of RNN steps due to parameter sharing across all RNN steps.}
For a relatively compact model with \XH{$|\mathbf{h}|=20$}, Fig. \ref{fig.1Dep_nsteps}(a) shows that the relative errors of training and testing samples decrease significantly when the number of RNN steps reaches 4, a critical number of RNN steps, \XH{which is expected since increasing the number of RNN steps enables the model to extract more accurate path-dependent features from longer-term stress-strain history.
As the number RNN steps further increases, the relative errors of training and testing samples remain at a plateau}, with around 0.4$\%$ and 1.6$\%$ relative errors, respectively.
The plateau indicates that further increasing the number of RNN steps does not improve the \XH{quality of the machine-inferred ISVs and thus the} model accuracy\XH{, which could be potentially limited by $|\hat{\mathbf{z}}|$ or $|\mathbf{h}|$}.
For a more complex model with \XH{$|\mathbf{h}|=30$}, Fig. \ref{fig.1Dep_nsteps}(b) shows a similar convergence behavior but the critical number of RNN steps increases to 5. 
This shows that the number of RNN steps play an important role in model accuracy and performance. 
Unnecessarily increasing the model complexity may lead to an increase in the number of RNN steps for achieving the same level of accuracy.

\subsubsection{Effects of Internal State Dimension}
\begin{figure}[htp]
\centering
    \begin{subfigure}{0.49\textwidth}
        \centering
        \includegraphics[width=1\linewidth]{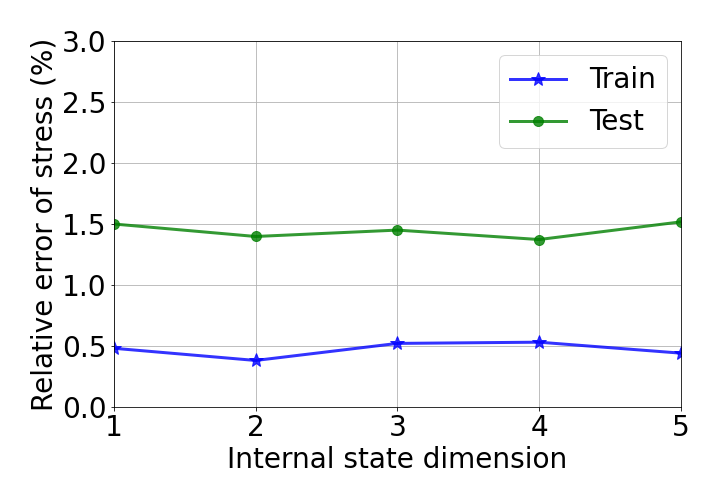}
        \vspace{-2\baselineskip}
        \caption{}
    \end{subfigure}
    \begin{subfigure}{0.49\textwidth}
        \centering
        \includegraphics[width=1\linewidth]{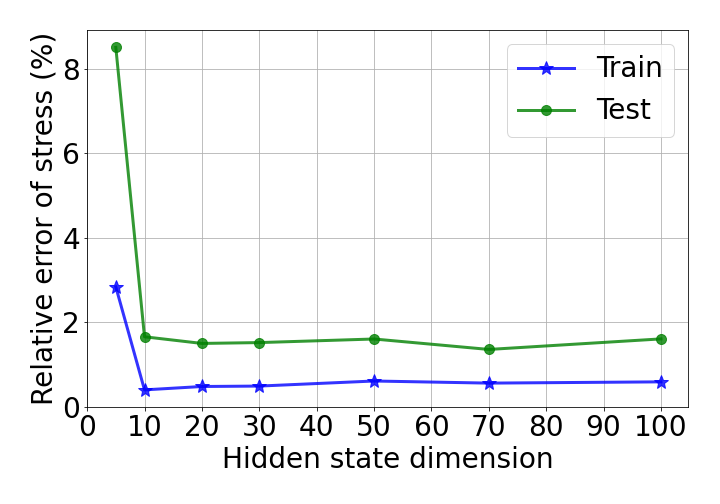}
        \vspace{-2\baselineskip}
        \caption{}
    \end{subfigure}
\caption{(a) Effects of the internal state dimension on the accuracy of models with 5 RNN steps and \XH{a hidden state dimension $|\mathbf{h}|=20$}. (b) Effects of the network complexity \XH{(hidden state dimension)} on the accuracy of models with a scalar internal state variable and 5 RNN steps.}\label{fig.1Dep_z_params}
\end{figure}

\XH{The internal state dimension ($|\hat{\mathbf{z}}|$) has a direct impact on the quality of the machine-inferred ISVs and thus the model performance. If $|\hat{\mathbf{z}}|$ is too small, the TCRNN model cannot capture all important thermodynamically consistent path-dependent features even if the number of RNN steps and model complexity are sufficient.} 
In this example, TCRNN models with 5 RNN steps, \XH{$|\mathbf{h}|=20$} and various internal state dimensions are examined. Fig. \ref{fig.1Dep_z_params}(a) shows that as the dimension of the ISV increases from 1 to 5, the relative errors of training and testing samples remain at a plateau, with around 0.5$\%$ and 1.5$\%$ relative errors, respectively. 
It indicates that a scalar ISV is sufficient for effectively capturing the path-dependent material behavior \XH{in this case}. 
The convergence of the model performance against the internal state dimension is particularly important. 
In practice, the internal state dimension of path-dependent materials is often unknown a priori. 
The convergence property shows that the TCRNN model remains accurate and robust even if an excessive internal state dimension is prescribed.
This convergence property also allows one to identify the optimal $|\hat{\mathbf{z}}|$ given measurable material states of path-dependent materials.

\subsubsection{Effects of Model Complexity}
\XH{The model complexity represented by the hidden state dimension is another important factor influencing the quality of the machine-inferred ISVs and model performance since the ISVs are inferred from the hidden states that directly capture the essential stress-strain path-dependent features. If the hidden state dimension is too small, important stress-strain path-dependent features could be lost,  leading to inaccurate machine-inferred ISVs and poor model performance.}
In this example, \XH{the effects the TCRNN model complexity (hidden state dimension) are investigated. The TCRNN models examined have a scalar ISV, 5 RNN steps, various hidden state dimensions ranging from 5 to 100.}
Fig. \ref{fig.1Dep_z_params}(b) shows that \XH{as the hidden state dimension increases,} the relative errors of training and testing samples decrease and then reach a plateau, with around 0.5$\%$ and 1.5$\%$ relative errors, respectively. 
This shows that a compact network is sufficient to achieve a satisfactory accuracy \XH{in this example} and increasing the model complexity does not improve the model accuracy.

\subsubsection{Model Generalization}
In the following examples, three variables are considered to evaluate the generalization performances of the TCRNN model, including the loading strain per cycle, the unloading strain per cycle, and the number of (loading-unloading) cycles. The TCRNN model with 15 RNN steps, a scalar ISV, and $|\mathbf{h}|=30$ is employed.

\begin{figure}[htp]
    \centering
    \includegraphics[width=1\linewidth]{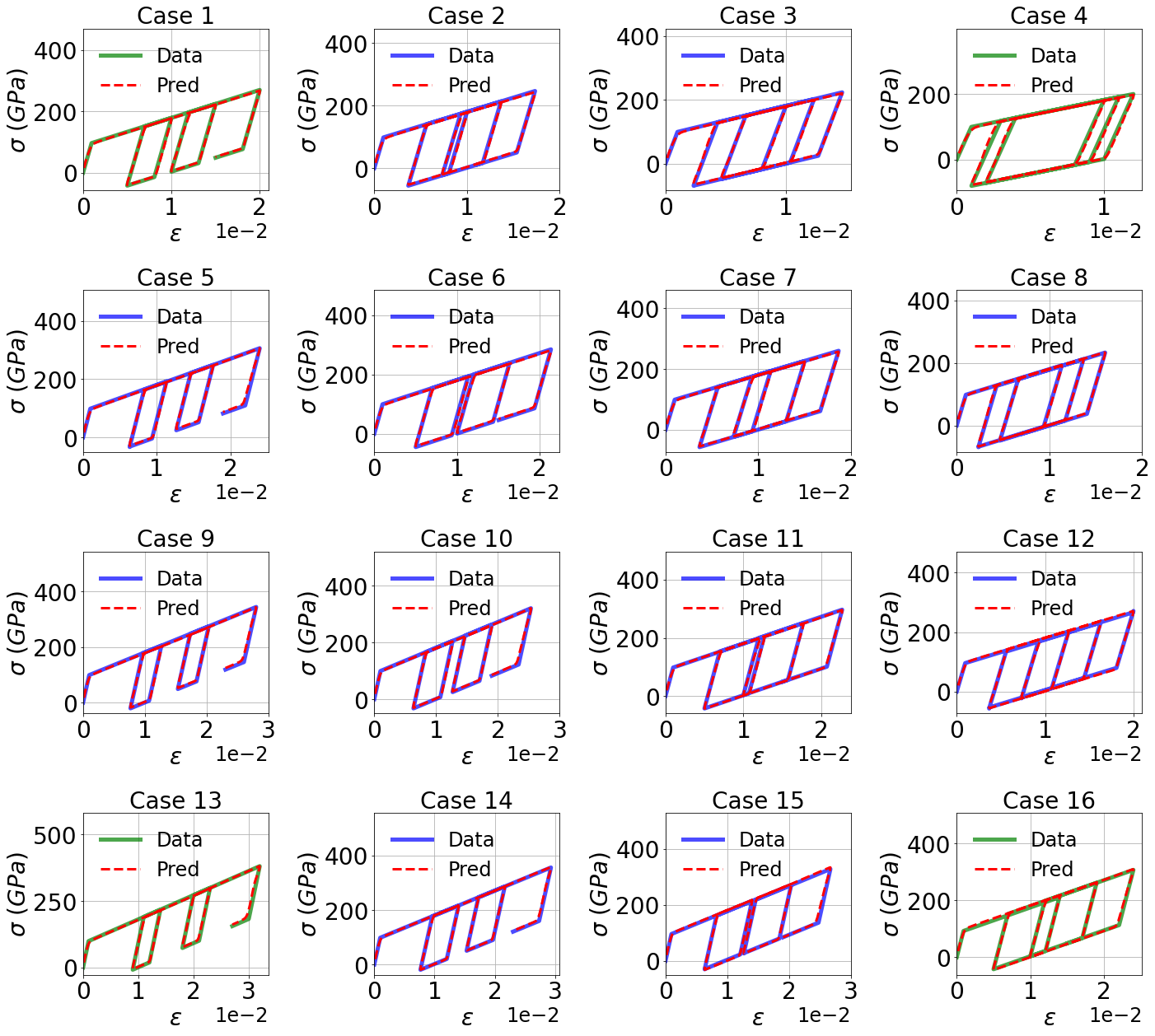}
    \caption{Comparison between the predicted stress of the TCRNN model and the data. The parameter space is constituted by the loading strain per cycle and the unloading strain per cycle. From top to bottom, the loading strain per cycle increases from $10^{-2}$ to $1.4 \times 10^{-2}$. From left to right, the unloading strain per cycle increases from $5 \times 10^{-3}$ to $5.5 \times 10^{-3}$. the The training data are denoted by the green color lines whereas the testing cases are denoted by the blue color lines. The predictions are denoted by the red dash lines.}\label{fig.1Dep_case1_s}
\end{figure}

In the first test, we consider a two-dimensional parameter space constituted by the loading strain per cycle and the unloading strain per cycle. 
The dataset contains 16 cases with the same number of loading-unloading cycles but with different loading and unloading strains per cycle. 
Fig. \ref{fig.1Dep_case1_s} shows the comparison between the predicted stress and the data, where case 1, 4, 9, and 16, located at the corners in the figure, are used for training with the data marked with the green solid lines, and the remaining cases are used for testing with the data marked with the blue solid lines. 
From top to bottom, the loading strain per cycle increases from $10^{-2}$ to $1.4 \times 10^{-2}$. 
From left to right, the unloading strain per cycle increases from $5 \times 10^{-3}$ to $5.5 \times 10^{-3}$. 
The mean relative errors of the training and testing cases are 3.1$\%$ and 2.8$\%$, respectively.
The good agreement between the data and the predictions demonstrates that the trained TCRNN model can successfully predict the testing cases within the prescribed parameter space.

\begin{figure}[htp]
    \centering
    \includegraphics[width=1\linewidth]{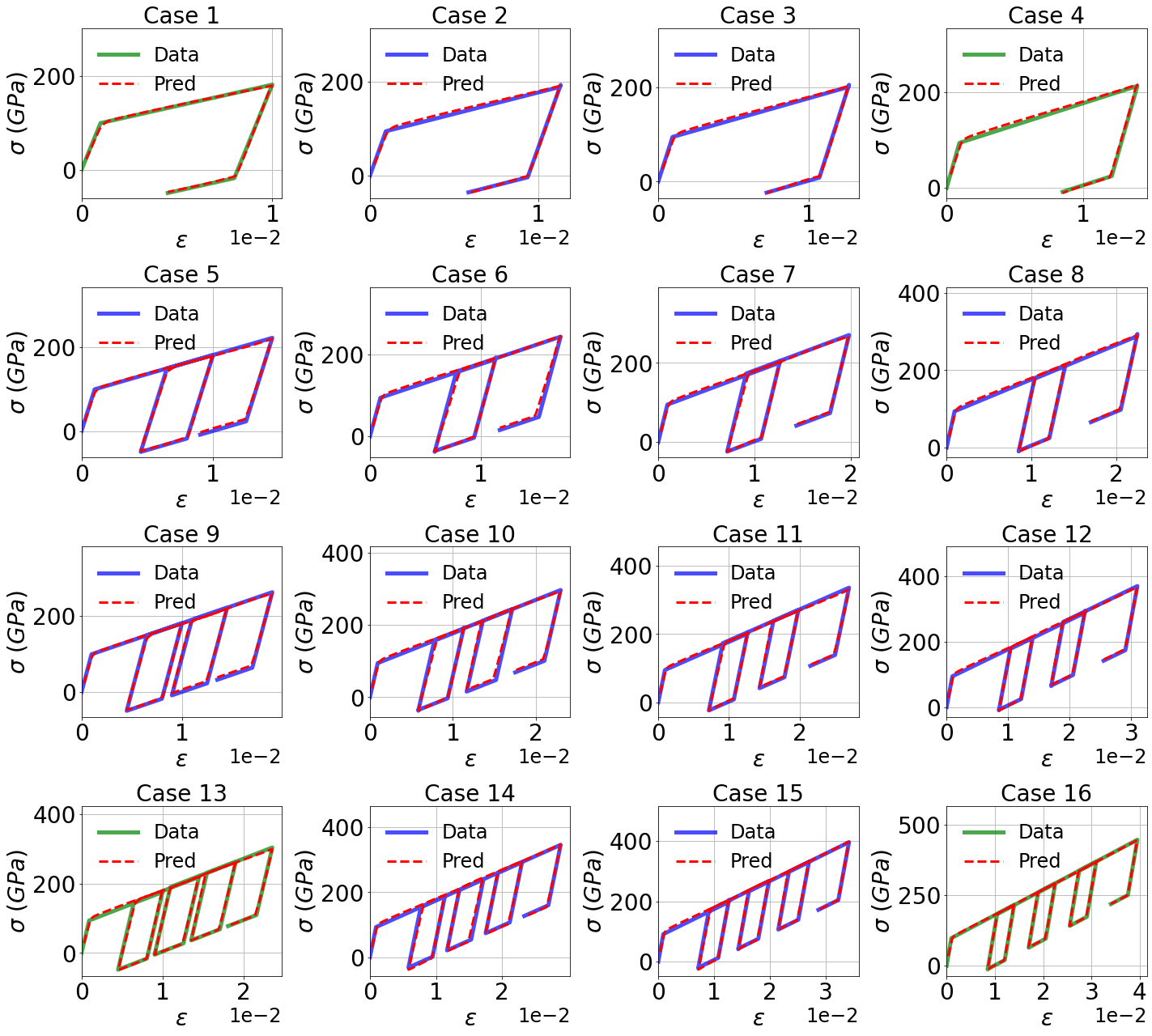}
    \caption{Comparison between the predicted stress of the TCRNN model and the data. The parameter space is constituted by the loading strain per cycle and the number of loading-unloading cycles. From top to bottom, the number of loading-unloading cycles increases from 1 to 4. From left to right, the loading strain per cycle increases from $10^{-2}$ to $1.4 \times 10^{-2}$. The training data are denoted by the green color lines whereas the testing cases are denoted by the blue color lines. The predictions are denoted by the red dash lines.}\label{fig.1Dep_case2_s}
\end{figure}

In the second test, we consider a two-dimensional parameter space constituted by the number of loading-unloading cycles and the loading strain per cycle. 
The dataset contains 16 cases with the same unloading strain per cycle, $5.5 \times 10^{-3}$.
Fig. \ref{fig.1Dep_case2_s} shows the comparison between the predicted stress and the data, where case 1, 4, 9, and 16, located at the corners in the figure, are used for training with the data marked with the green solid lines, and the remaining cases are used for testing with the data marked with the blue solid lines. 
From top to bottom, the number of loading-unloading cycles increases from 1 to 4. 
From left to right, the loading strain per cycle increases from $10^{-2}$ to $1.4 \times 10^{-2}$.
The mean relative errors of the training and testing cases are 2.3$\%$ and 3.4$\%$, respectively.
The good agreement between the data and the predictions further demonstrates the strong generalization ability of the TCRNN constitutive model.

\subsection{Modeling Soil under Cyclic Shear Loading}\label{sec:result_soil}
\begin{figure}[htp]
\centering
    \begin{subfigure}{1\textwidth}
        \centering
        \includegraphics[width=1\linewidth]{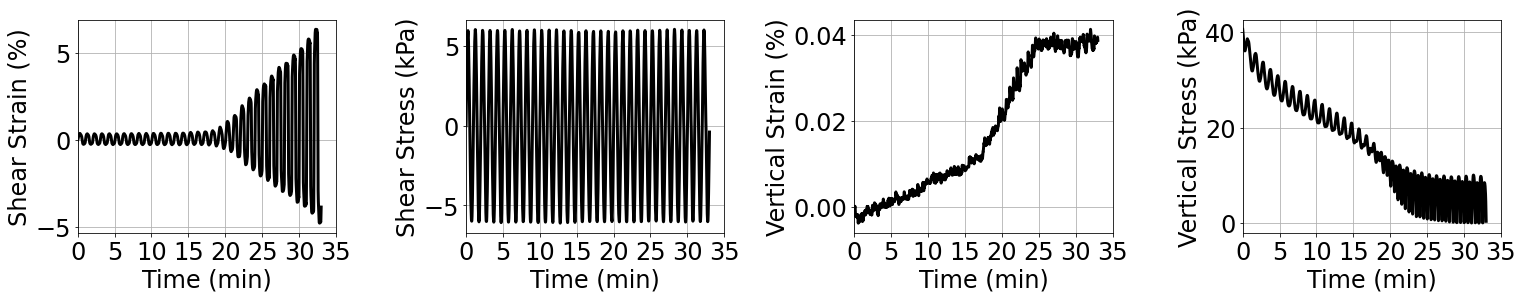}
        \vspace{-1.5\baselineskip}
        \caption{CSR=0.15}
    \end{subfigure}
    \begin{subfigure}{1\textwidth}
        \centering
        \includegraphics[width=1\linewidth]{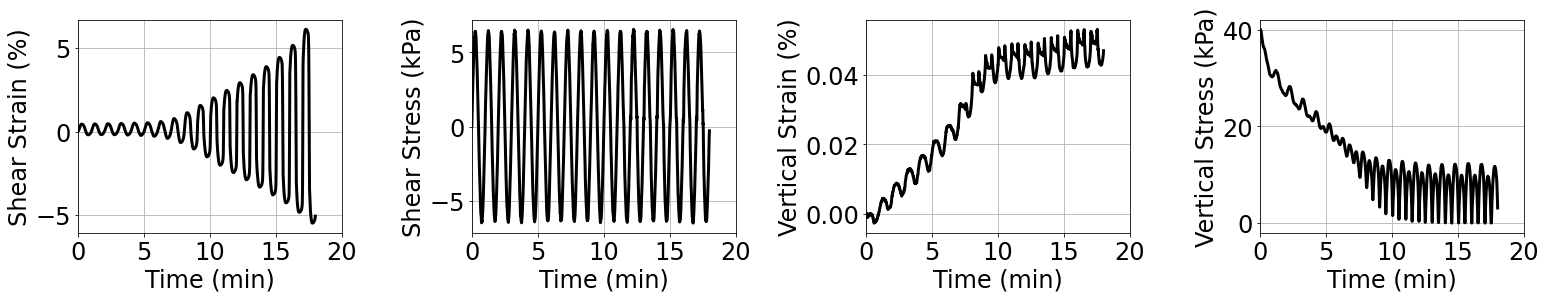}
        \vspace{-1.5\baselineskip}
        \caption{CSR=0.16}
    \end{subfigure}
    \begin{subfigure}{1\textwidth}
        \centering
        \includegraphics[width=1\linewidth]{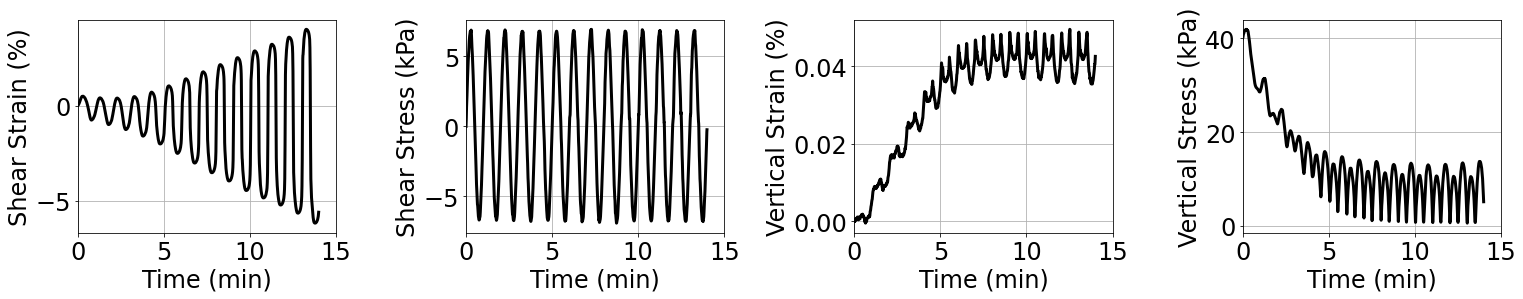}
        \vspace{-1.5\baselineskip}
        \caption{CSR=0.17}
    \end{subfigure}
\caption{Experimental data of undrained soil under cyclic simple shear loading with: (a) CSR=0.15; (b) CSR=0.16; (c) CSR=0.17.}\label{fig.soil_data}
\end{figure}

The effectiveness of the proposed TCRNN constitutive model is further evaluated by modeling undrained soil under cyclic shear loading \cite{bastidas2016ottawa,ghoraiby2020physical}. 
The experimental data is collected from the undrained soil samples under initial triaxial confinement of $40 kPa$ and cyclic shear loading. A cyclic stress ratio (CSR) is defined as the ratio of the maximum shear stress to the initial vertical stress. 
The experimental data contains the shear strain, the vertical strain, the shear stress, and the vertical stress. 
Fig. \ref{fig.soil_data} shows the experimental data with a CSR of 0.15, 0.16, and 0.17. 
The stress-strain relationships are highly nonlinear and path-dependent due to coupling effects of changes in volume, matric suction, degree of saturation, effective stress, shear modulus, etc. \cite{rong2021undrained}. 
Modeling such path-dependent material behaviors by a phenomenological approach is challenging and complicated, which often relies on certain phenomenological ISVs.

Given only the stress-strain data, data-driven models \XH{and phenomenological models} that require pre-defined ISVs cannot be applied. 
In contrast, the proposed TCRNN model can be effectively applied since it only requires measurable material states, and the model is capable of inferring essential ISVs from the measurable material states by following the thermodynamics principles. 

The TCRNN \XH{based on the increment of ISVs ($\Delta \hat{\mathbf{z}}$)}, as shown in Fig. \ref{fig.tcrnn_dz}(b), is employed in this example.
A GRU is used to infer the ISV $\hat{\mathbf{z}}$ in Eq. \eqref{eq.tcrnn_rnn_bar} and describe its evolution by following the thermodynamics second law in Eq. \eqref{eq.2nd_law}.
\XH{The GRU consists of one hidden layer for all affine transformations in Eq. \eqref{eq.rnn-gru} with the model complexity represented by the hidden state dimension, $|\mathbf{h}|$.}
Since the training data contains only stresses and strains, the loss function in Eq. \eqref{eq.tcrnn_loss3} is employed with $\beta_1=\beta_2=1$. 
The experimental data with a CSR of 0.15 and 0.17 are used for training, while the experimental data with a CSR of 0.16 is used for testing. 
The effects of the number of RNN steps, the internal state dimension, and the model complexity on the model performance are investigated. 

\begin{figure}[htp]
\centering
    \begin{subfigure}{0.45\textwidth}
        \centering
        \includegraphics[width=1\linewidth]{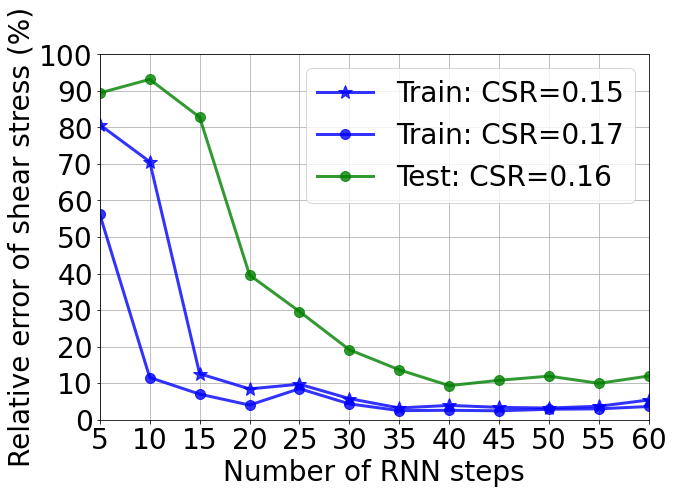}
        \vspace{-1.5\baselineskip}
        \caption{}
    \end{subfigure}
    \begin{subfigure}{0.45\textwidth}
        \centering
        \includegraphics[width=1\linewidth]{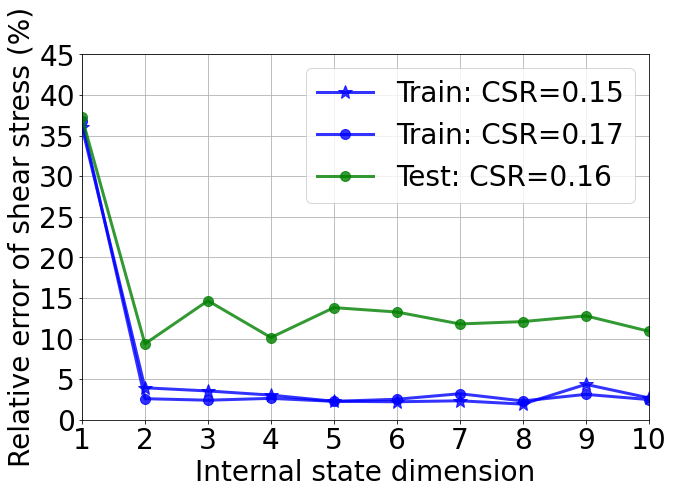}
        \vspace{-1.5\baselineskip}
        \caption{}
    \end{subfigure}
    \begin{subfigure}{0.45\textwidth}
        \centering
        \includegraphics[width=1\linewidth]{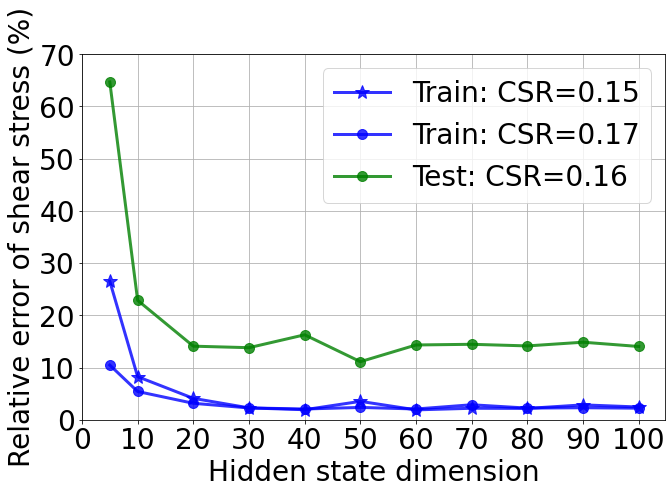}
        \vspace{-1.5\baselineskip}
        \caption{}
    \end{subfigure}
\caption{Effects of various parameters on model accuracy: (a) the number of RNN steps (with an internal state dimension $|\hat{\mathbf{z}}|=2$ and \XH{a hidden state dimension $|\mathbf{h}|=30$}); (b) the internal state dimension (with 40 RNN steps and \XH{a hidden state dimension $|\mathbf{h}|=30$}); (c) the model complexity (with an internal state dimension $|\hat{\mathbf{z}}|=5$ and 40 RNN steps).}\label{fig.soil_parameter_effects}
\end{figure}

Given TCRNN models with an internal state dimension \XH{($|\hat{\mathbf{z}}|$)} of 2 and \XH{a hidden state dimension ($|\mathbf{h}|$) of 30}, the number of RNN steps is varied from 5 to 60 and its influences on the model prediction accuracy are shown in Fig. \ref{fig.soil_parameter_effects}(a)
As the number of RNN steps increases, the relative errors of training and testing samples decrease and eventually converge to a plateau, with values around 3$\%$ and 11$\%$, respectively.
The plateau indicates that further increasing the number of RNN steps does not improve the model accuracy.

The internal state dimension \XH{($|\hat{\mathbf{z}}|$)} required to effectively model the path-dependent material behaviors is unknown a priori, which depends on the complexity of the path-dependent behaviors. 
Here, we investigate the effects of $|\hat{\mathbf{z}}|$ on model prediction accuracy, which is varied from 1 to 10, while the number of RNN steps and \XH{$|\mathbf{h}|$} are fixed as 40 and 30, respectively. 
Fig. \ref{fig.soil_parameter_effects}(b) shows that the relative errors of training and testing samples are large when the machine-inferred ISV is a scalar, indicating that a scalar ISV is insufficient to capture all essential path-dependent features.
As $|\hat{\mathbf{z}}|$ increases, the relative errors of training and testing samples decrease and then reach a plateau, with values around 2.7$\%$ and 12$\%$, respectively, which shows that the TCRNN model remains accurate and robust even if an excessive $|\hat{\mathbf{z}}|$ is prescribed.
The convergence behavior also allows one to identify the optimal $|\hat{\mathbf{z}}|$, around 2 in this example, given the measurable material states of path-dependent materials.


Lastly, $|\hat{\mathbf{z}}|$ and the number of RNN steps are fixed as 5 and 40, respectively, while \XH{$|\mathbf{h}|$ is varied from 5 to 100 to investigate the effects of model complexity ($|\mathbf{h}|$) on model performance}.
Fig. \ref{fig.soil_parameter_effects}(c) shows that the relative errors of training and testing samples decrease as \XH{$|\mathbf{h}|$} increases and eventually reach a plateau, with values around 2.8$\%$ and 14$\%$, respectively. 
The plateaus in the convergence curves indicate that further increasing the model complexity does not improve the model accuracy.

\begin{figure}[htp]
\centering
    \begin{subfigure}{0.325\textwidth}
        \centering
        \includegraphics[width=1\linewidth]{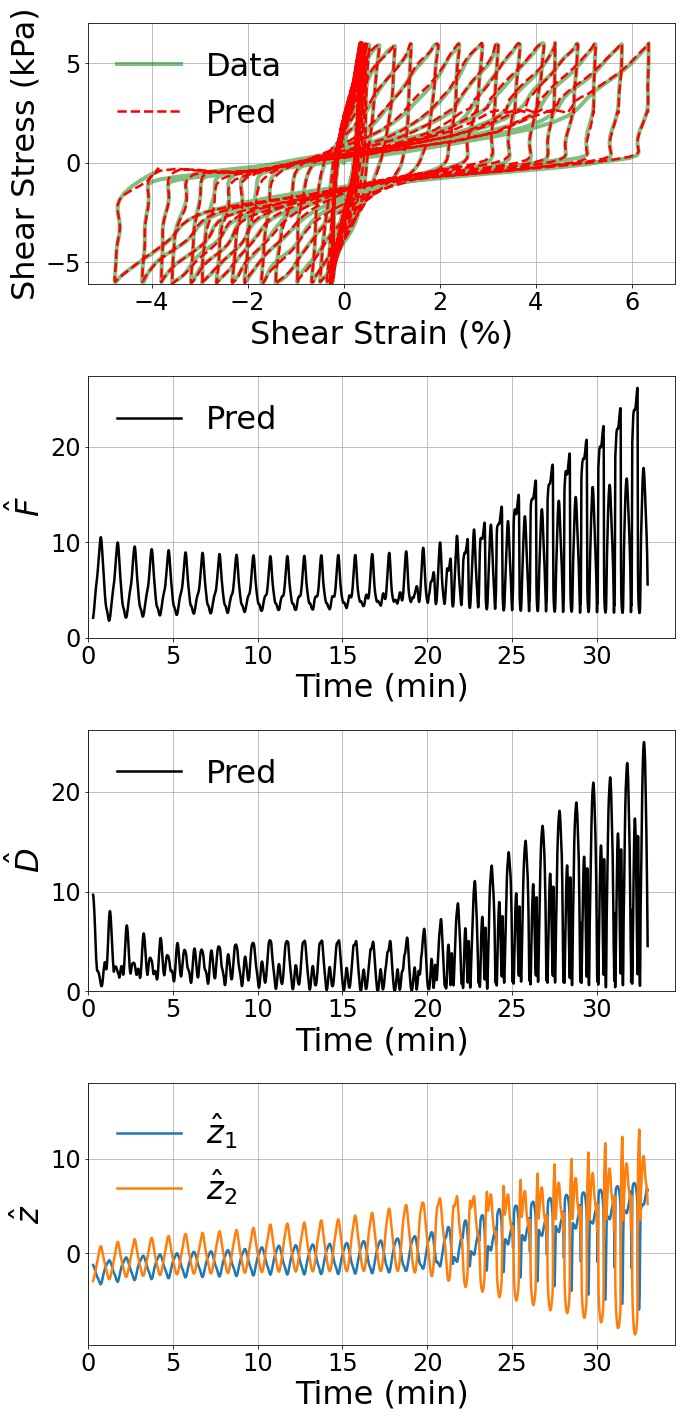}
        \caption{Training: CSR=0.15}
    \end{subfigure}
    \begin{subfigure}{0.325\textwidth}
        \centering
        \includegraphics[width=1\linewidth]{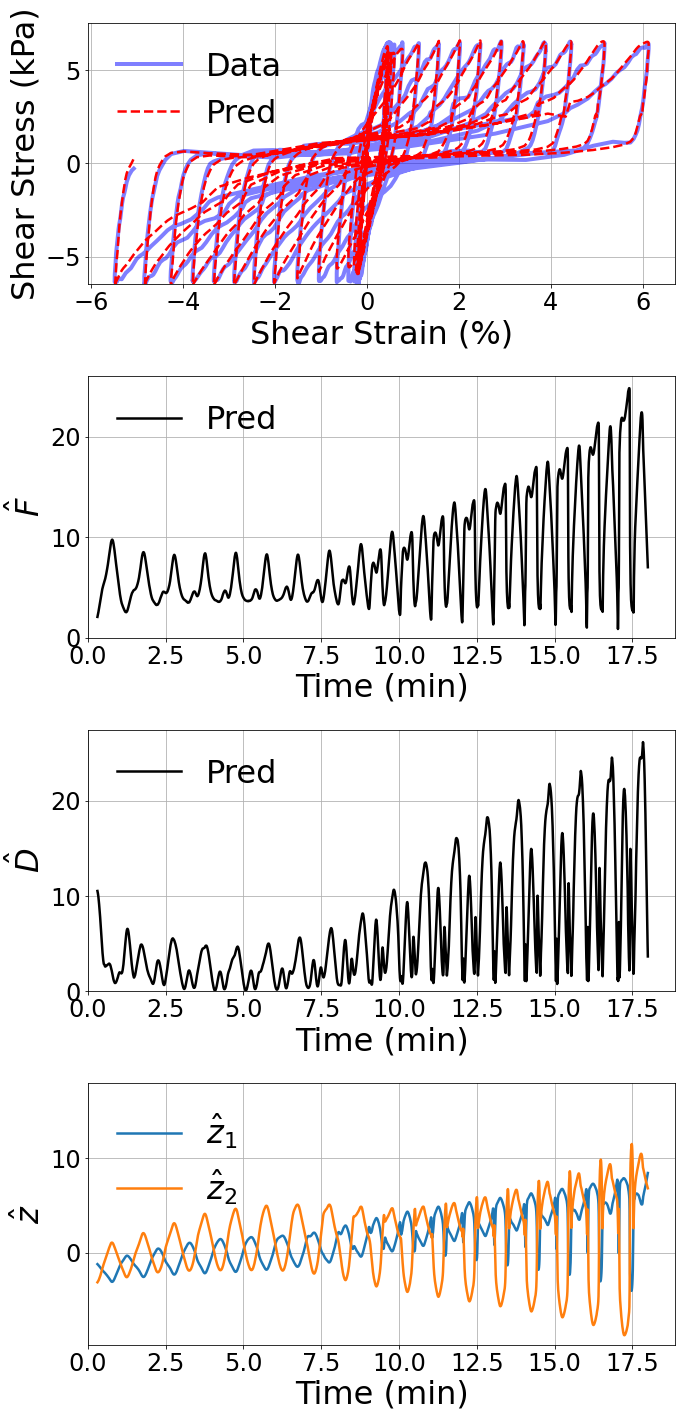}
        \caption{Testing: CSR=0.16}
    \end{subfigure}
    \begin{subfigure}{0.325\textwidth}
        \centering
        \includegraphics[width=1\linewidth]{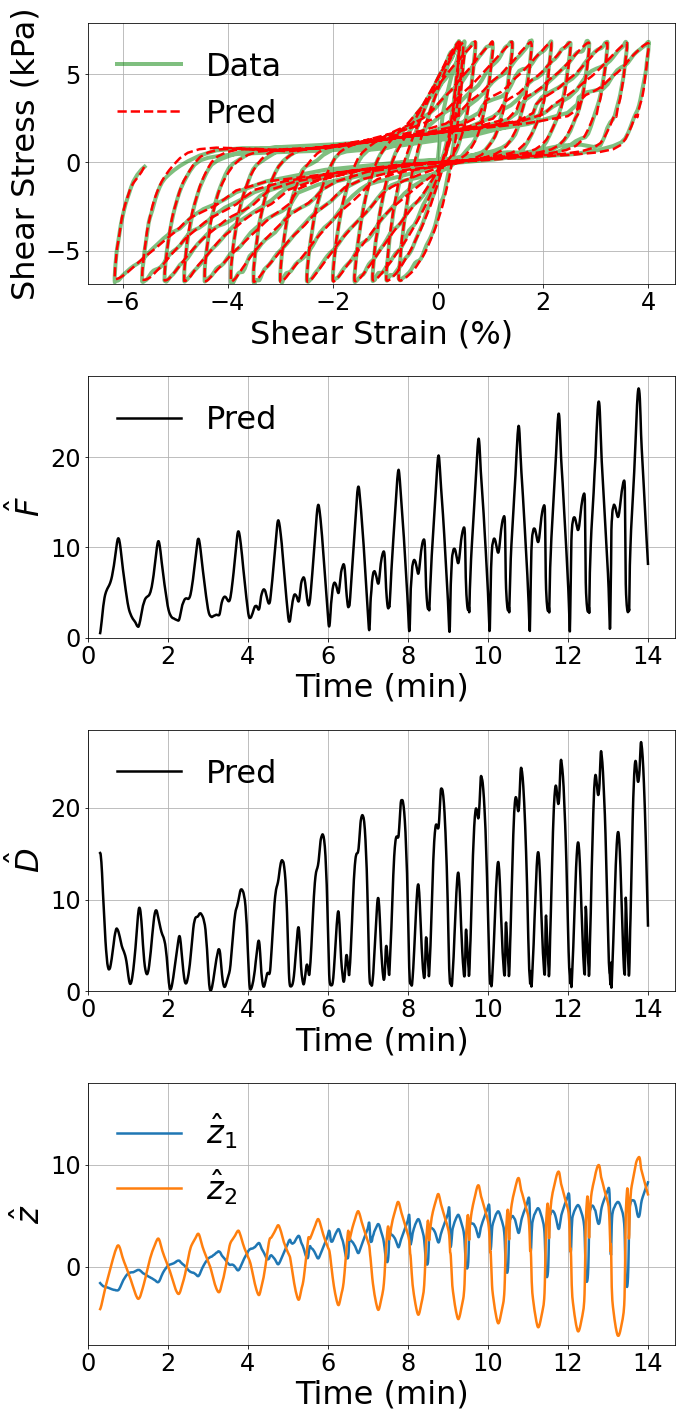}
        \caption{Training: CSR=0.17}
    \end{subfigure}
\caption{Comparison of predictions of the TCRNN with data: (a) the training case with a CSR=0.15; (b) the testing case with a CSR=0.16; (c) the training case with a CSR=0.17. The first row compares shear stress-strain relationships. The second row shows the predicted Helmholtz free energy. The third row shows the predicted dissipation rate. The last row shows the machine-learned internal state variables. The TCRNN model employed has 40 RNN steps, an internal state dimension $|\hat{\mathbf{z}}|=2$, and \XH{a hidden state dimension $|\mathbf{h}|=30$}.}\label{fig.soil_train_test}
\end{figure}

Fig. \ref{fig.soil_train_test} compares shear stress experimental data with the predictions of the trained TCRNN model that employs 40 RNN steps, $|\hat{\mathbf{z}}|=2$, and \XH{$|\mathbf{h}|=30$}. 
The relative errors of training and testing samples are around 3.2$\%$ and 9.4$\%$, respectively. 
It shows that the TCRNN model is able to learn the path-dependent material behaviors from the measurable material states under given loading conditions and effectively predicts the path-dependent responses under untrained loading conditions, further demonstrating the generalization ability and effectiveness of the TCRNN model in practical applications. 
Further, the trained TCRNN material model is thermodynamically consistent, which is verified by the non-negative predicted free energy and the predicted dissipation rates that satisfy the thermodynamics second law, as shown in the second and the third rows of Fig. \ref{fig.soil_train_test}, respectively.
The histories of the machine-learned ISVs are shown in the last row of Fig. \ref{fig.soil_train_test}, revealing interesting path-dependent patterns similar to the behaviors of the predicted free energy and dissipation rate.

\section{Conclusions}\label{sec:conclusion}
In this study, we introduced a machine-learned internal state variable (ISV) approach for data-driven modeling of path-dependent materials, which is thermodynamically consistent and relies purely on the measurable material states. 
The proposed TCRNN constitutive models consist of two main components: an RNN that infers ISVs (Eq. \eqref{eq.tcrnn_rnn_bar}) and describes their evolution by following the thermodynamics second law (Eq. \eqref{eq.2nd_law}), and a DNN that predicts the Helmholtz free energy (Eq. \eqref{eq.tcrnn_nn_bar}) given strain, ISVs, and temperature (for non-isothermal processes). 
Two TCRNN constitutive models are developed, one based on \XH{the time rates of ISVs ($\dot{\hat{\mathbf{z}}}$)}, as shown in Fig. \ref{fig.tcrnn_dzdt}, and the other one based on \XH{the increments of ISVs ($\Delta \hat{\mathbf{z}}$)}, as shown in Fig. \ref{fig.tcrnn_dz}. 
The latter model shows an enhanced efficiency as it utilizes an approximation of $\dot{\hat{\mathbf{z}}}$ for the calculation of dissipation rate and avoids time-consuming differentiation of the RNN outputs with respect to all RNN inputs.
Model robustness and accuracy is enhanced by introducing \textit{stochasticity} to the training data to account for uncertainties of input conditions in the testing.

In the demonstration of modeling elasto-plastic materials, the parametric study shows that the model accuracy converges as the number of RNN steps, the internal state dimension, and the model complexity increase. 
\XH{All these factors} play an important role in the model performance. 
Given path-dependent material behaviors, there exists an \XH{optimal} internal state dimension to capture the essential path-dependent features by the TCRNN model. 
It has been shown that the TCRNN model remains accurate and robust even if an excessive internal state dimension is prescribed.
The monotonic correlation between the machine-inferred and the phenomenological ISV \XH{of the elasto-plastic material} demonstrates that the TCRNN constitutive model can infer mechanistically and thermodynamically consistent ISVs. 
The proposed TCRNN constitutive model is shown to remain robust against various strain increments and have strong generalization capabilities.

The effectiveness of the proposed TCRNN constitutive model is further demonstrated by modeling undrained soil under cyclic shear loading using experimental data, where only measurable material states (stresses and strains) are available. 
A similar convergence behaviors of the model accuracy are observed from a parametric study of the number of RNN steps, the internal state dimension, and the model complexity. 
The generalization capability of the TCRNN constitutive model is demonstrated by the effective prediction of the thermodynamically consistent response of undrained soil under the loading conditions different from the ones used in training, which reveals the promising potential of the proposed method to model complex path-dependent materials behaviors in real applications.

The proposed TCRNN constitutive model is general and applicable to model a wide range of path-dependent materials. 
It is efficient and can be applied to accelerate large-scale multi-scale simulations with complex microstructures and path-dependent material systems. 
To investigate reliability of model predictions, a future extension would be to integrate uncertainty quantification \cite{smith2013uncertainty} into the proposed TCRNN model.

\section*{Acknowledgements}
The support of this work by the DOE Nuclear Engineering University Program under Award Number DE-NE0008951 to the University of California, San Diego is very much appreciated.


\bibliography{reference}

\end{document}